\input gtmacros
\input amsnames
\input amstex

\catcode`\@=12        

\input gtmonout
\pagenumbers{507}{553}
\papernumber{21}
\received{9 November 1998}\revised{27 January 1999}
\published{21 October 1999}

\let\\\par
\def\topmatter{\relax}

\let\gttitle\title
\def\title#1\endtitle{\gttitle{#1}}
\let\gtauthor\author
\def\author#1\endauthor{\gtauthor{#1}}
\let\gtaddress\address
\def\address#1\endaddress{\gtaddress{#1}}
\let\gtemail\email
\def\email#1\endemail{\gtemail{#1}}
\def\subjclass#1\endsubjclass{\primaryclass{#1}}
\let\gtkeywords\keywords
\def\keywords#1\endkeywords{\gtkeywords{#1}}
\def\heading#1\endheading{{\def\S##1{\relax}\def\\{\relax\ignorespaces}
    \section{#1}}}
\def\head#1\endhead{\heading#1\endheading}

\def\subhead#1\endsubhead{\sh{#1}}
\def\subsubhead#1\endsubsubhead{\sh{#1}}
\def\specialhead#1\endspecialhead{\sh{#1}}
\def\demo#1{\rk{#1}\ignorespaces}
\def\enddemo{\ppar}

\let\oldsq\sq
\def\qed{\ifmmode\quad\oldsq\else\hbox{}\hfill$\oldsq$\par\goodbreak\rm\fi}  
\def\proclaim#1{\rk{#1}\sl\ignorespaces}
\def\endproclaim{\rm\ppar}
\def\cite#1{[#1]}
\newcount\itemnumber
\def\roster{\items\itemnumber=1}
\def\endroster{\enditems}
\let\itemold\item
\def\item{\itemold{{\rm(\number\itemnumber)}}%
\global\advance\itemnumber by 1\ignorespaces}
\def\S{section~\ignorespaces}  
\def\date#1\enddate{\relax}
\def\thanks#1\endthanks{\relax}   
\def\dedicatory#1\enddedicatory{\relax}  
\let\footnote\plainfootnote

\def\Refs{\ppar{\large\bf References}\ppar\bgroup\leftskip=25pt
\frenchspacing\parskip=3pt plus2pt\small}       
\def\endRefs{\egroup}
\def\widestnumber#1#2{\relax}
\def\endrefitem{}
\def\refdef#1#2#3{\def#1{\leavevmode\unskip\endrefitem#2\def\endrefitem{#3}}}
\def\ref{\par}
\def\endref{\endrefitem\par\def\endrefitem{}}
\refdef\key{\noindent\llap\bgroup[}{]\ \ \egroup}
\refdef\no{\noindent\llap\bgroup[}{]\ \ \egroup}
\refdef\by{\bf}{\rm, }
\refdef\manyby{\bf}{\rm, }
\refdef\paper{\it}{\rm, }
\refdef\book{\it}{\rm, }
\refdef\jour{}{ }
\refdef\vol{}{ }
\refdef\yr{$(}{)$ }
\refdef\ed{(}{ Editor) }
\refdef\publ{}{ }
\refdef\inbook{from: ``}{'', }
\refdef\pages{}{ }
\refdef\page{}{ }
\refdef\paperinfo{}{ }
\refdef\bookinfo{}{ }
\refdef\publaddr{}{ }
\refdef\eds{(}{ Editors)}
\refdef\bysame{\hbox to 3 em{\hrulefill}\thinspace,}{ }
\refdef\toappear{(to appear)}{ }
\refdef\issue{no.\ }{ }
\newcount\refnumber\refnumber=1
\def\refkey#1{\expandafter\xdef\csname cite#1\endcsname{\number\refnumber}%
\global\advance\refnumber by 1}
\def\cite#1{[\csname cite#1\endcsname]}
\def\Cite#1{\csname cite#1\endcsname}  
\def\key#1{\noindent\llap{[\csname cite#1\endcsname]\ \ }}

\refkey {BFSV}
\refkey {B}
\refkey {BQ}
\refkey {CFS}
\refkey {DW}
\refkey {F}
\refkey {FQ}
\refkey {FK}
\refkey {GK}
\refkey {Ka}
\refkey {Ke}
\refkey {K}
\refkey {Ku}
\refkey {L}
\refkey {McL}
\refkey {M}
\refkey {Q}
\refkey {QTP}
\refkey {RT}
\refkey {Sa}
\refkey {S}
\refkey {T}
\refkey {W}
\refkey {Wh}
\refkey {Y}

\def\hom{\text{\rm hom}}
\def\units{\text{\rm units}}
\def\sq{\diamond}
\def\id{\text{\rm id}}
\def\pt{\text{\rm pt}}
\def\br#1{\bar{#1}}

\topmatter
\title Group categories and their field theories\endtitle
\author Frank Quinn\endauthor
\abstract
A group--category is an additively semisimple category with a monoidal
product structure in which the simple objects
are invertible. For example in the category of representations of a group,
1--dimensional representations are the
invertible simple objects. This paper gives a detailed exploration of
``topological quantum field theories'' for
group--categories, in hopes of finding clues to a better understanding of
the general situation. Group--categories are
classified in several ways extending results of Fr\"olich and Kerler.
Topological field theories based on homology and
cohomology are constructed, and these are shown to include theories
obtained from group--categories by
Reshetikhin--Turaev constructions. Braided--commutative categories most
naturally give theories on 4--manifold
thickenings of 2--complexes; the usual 3--manifold theories are obtained from
these by normalizing them (using results
of Kirby) to depend mostly on the boundary of the thickening. This is
worked out for group--categories, and in
particular we determine when the normalization is possible and when it is not.

\endabstract

\asciiabstract{A group-category is an additively semisimple category
with a monoidal product structure in which the simple objects are
invertible. For example in the category of representations of a group,
1-dimensional representations are the invertible simple objects. This
paper gives a detailed exploration of "topological quantum field
theories" for group-categories, in hopes of finding clues to a better
understanding of the general situation. Group-categories are
classified in several ways extending results of Froelich and Kerler.
Topological field theories based on homology and cohomology are
constructed, and these are shown to include theories obtained from
group-categories by Reshetikhin-Turaev constructions.
Braided-commutative categories most naturally give theories on
4-manifold thickenings of 2-complexes; the usual 3-manifold theories
are obtained from these by normalizing them (using results of Kirby)
to depend mostly on the boundary of the thickening. This is worked out
for group-categories, and in particular we determine when the
normalization is possible and when it is not.}

\address Department of Mathematics, Virginia Tech, Blacksburg VA
24061-0123, USA\endaddress
\email quinn@math.vt.edu\endemail
\date November 1998\enddate

\primaryclass{18D10}\secondaryclass{81R50, 55B20}

\keywords Topological quantum field theory, braided category
\endkeywords

\maketitle

\cl{\small\it 
Dedicated to Rob Kirby, on the occasion of his $60^{th}$ birthday}

\catcode`\@=\active 

\head Introduction\endhead
There is a close connection between monoidal categories and
low-dimensional modular topological field theories.  Specifically,
symmetric monoidal categories correspond to field theories on
2--dimensional CW complexes [\Cite{B}, \Cite{Q}]; monoidal categories
correspond to theories on 3--manifolds with boundary, and tortile
(braided--commutative) categories correspond to theories on
4--dimensional\break thickenings of 2--complexes. These last can
usually be normalized to give theories on extended 3--manifolds, and
this is the most familiar context [\Cite{RT}, \Cite{T}, \Cite{Ke},
\Cite{Sa}, \Cite{W}].  Particularly interesting braided categories are
obtained from representations of ``quantum groups'' at roots of unity,
cf [\Cite{L}, \Cite{Ka}], and analogous symmetric mod $p$ categories
were defined by Gelfand and Kazhdan \cite{GK}.

This subject has produced a voluminous literature but not a lot of new
information. Presumably we do not yet
understand the geometric significance, wider contexts, methods of
computation, etc, well enough to effectively
exploit these theories.  This paper presents a class of
examples in which everything can be worked out in detail, as a source  of
clues for the general case.
Descriptions of the categories gives a  connection to recent work on
classifying spaces. The field theories
turn out to be special cases of constructions using homology of CW
complexes, or more generally cohomology of manifold
thickenings of CW complexes. This clarifies the nature of the objects on
which the fields are defined, and hints at
higher-dimensional versions. The examples illuminate the normalization
procedure used to pass to fields on
extended 3--manifolds. Finally group--categories occur as tensor factors of
the ``quantum''
categories (2.2.4), so understanding them is an essential ingredient of the
general case.

Finite
groups provide another class of examples that have been worked out in detail
[\Cite{FQ}, \Cite{Q}, \Cite{Y}], but these have not been so helpful.  Representations of
the group give a (symmetric) monoidal
category, and  a field theory (on all finite CW complexes) defined in terms of
homomorphisms of fundamental groups into the finite group. The restriction
of the field theory to 2--complexes is
the field theory corresponding to the representation category. However the
restriction of the field theory to
3--manifolds corresponds to the {\it double\/} of the category \cite{M}, not
the category itself. Constructions using
a double
 are much easier but also much less informative than the general case, so
this is a defect in this model.

 A {\it group category\/} is a
semisimple additive category with a product structure in which the
simple objects are invertible.  Isomorphism classes of
simple objects then form a group, called the ``underlying group'' of
the category. Section 2 begins with a slightly more precise definition
(2.1) and some examples. The conjectural
appearance of group--categories as tensor factors of quantum categories
(2.2.4) is particularly curious. Three views of
the classification of group--categories are then presented. The first and
only novel view (2.3) uses recent work on
classifying spaces of braided categories \cite{F} to give a
characterization in terms of spaces with two nonvanishing
homotopy groups. Specifically, group--categories over a ring $R$ with
underlying group $G$ correspond to spaces $E$ with
$\pi_d(E)=G$ and
$\pi_{d+1}(E)=\units(R)$. The cases $d=1, 2$, and $d\geq 3$ correspond to
monoidal, braided--commutative, and symmetric
categories respectively. The Postnikov decomposition  gives an equivalence
of this to
$k$--invariants in group cohomology. The second approach (2.4) derives a
category structure directly from  group
cohomology using cellular cochains in a model for the classifying space.
This approach was developed by Fr\"olich and
Kerler
\cite{FK}. The third approach (2.5) gives a ``numerical presentation'' for
the category. This is a format developed
for machine computation
[\Cite{BQ}, \Cite{QTP}], but in this case it  gives an explicit and efficient
low-level description.

Group cohomology in the context of topological field theories first
appeared in Dijkgraaf--Witten \cite{DW} as
lagrangians for fields with finite gauge group.  Their lagrangians lie in
$H^3(B_G)$, which we now see as classifying
monoidal (no commutativity conditions) group--categories. The field theory
they construct corresponds to the double of
the category.

Topological field theories based on homology with coefficients in a finite
group are
studied in section 3.  Suppose $G$ is a finite abelian group and $R$ a
ring.  State spaces of the $H_n$ theory are the free modules
$R[H_{n}(Y;G)]$. Induced
homomorphisms are defined by summing over $H_{n+1}$: if $X\supset Y_{1},
Y_{2}$ and $y\in H_{1}(Y_1;G)$ then
$$Z_{X}(y)=\Sigma_{\{x\in H_{n+1}(X;G)\mid\partial_{1}x=-y\}}\partial_{2}x.$$
We determine (3.1.3) exactly when this satisfies various field theory
axioms. The $H_1$ theory is the one that
connects with categories: on 2--complexes it corresponds to the standard
(untwisted) group--category. On 3--manifolds it
provides examples of field theories that are not modular. This illustrates
the role of doubling or extended
structures  in obtaining modularity on 3--manifolds. The higher-dimensional
versions are new, and suggest interesting
connections with classical algebraic topology.

Probably the eventual proper setting for field theories will be covariant
(homological),
but the current constructions are too rigid. In section 4 we restrict to
manifolds and consider the dual
cohomology-based theories. Here we can build in a twisting by  evaluating
group cohomology classes on
fundamental classes. Again we get examples for any $n$, and it is
the $n=1$ cases that relate to group--categories. Again homological
calculations determine when these satisfy field
axioms. For $n=1$ state spaces are associated to manifold with the homotopy
type of 1--complexes (we refer to these as
``thickenings'' of 1--complexes); induced homomorphisms come from
thickenings of 2--complexes, and corners used in
modular structures are thickenings of 0--complexes. The dimensions of these
thickenings depend on the type of category.
To establish notation we relate both fields and categories to spaces with
two homotopy groups. Let $E$ have
$\pi_d(E)=G$ and $\pi_{d+1}(E)=\units(R)$. Then $E$ determines a category
and a field theory:
\medskip
\hfil\vbox{\offinterlineskip
\halign{\vrule#&\strut\quad\hfil#\quad&\vrule#&\quad\hfil#\hfil
\quad&\vrule#&\quad\hfil#\hfil\quad&\vrule#\cr &$d$&&category structure
&&fields on & \cr
\noalign{\hrule}
&1&&associative &&$(3,2,1)$--thickenings &\cr
&2&&braided--commutative &&$(4,3,2)$--thickenings &\cr
&$\geq 3$&&symmetric &&$(d+2,d+1,d)$--thickenings &\cr
\noalign{\hrule}\cr
}}\hfil
\medskip
We show (4.3) that the field theory is in fact the one obtained by a
Reshetikhin--Turaev construction from the category.

Section 5 concerns field theories on 3--manifolds. The basic plan [\Cite{W},
\Cite{T}] is to start with a theory on
4--dimensional thickenings of 2--complexes, associated to a
braided--commutative category, and try to extract a theory
that depends only on the boundary of the thickening. The geometric
ingredient is the basis of the Kirby calculus
\cite{K}: a 3--manifold bounds a simply-connected 4--manifold, and this
4--manifold is well-defined up to connected sums
with
$CP^2$ and $\overline{CP}^2$. If we specify the index of the 4--manifold
then it is well-defined up to sums with
$CP^2\#\overline{CP}^2$. These connected sums change the induced
homomorphisms by multiplication by an element in $R$.
If the element associated to $CP^2\#\overline{CP}^2$ has an inverse square
root then we can use it to nomalize the
theory (tensor with an Euler characteristic theory) to be insensitive to
such sums. This gives a theory defined on
``extended'' 3--manifolds: manifolds together with an integer specifying the
index of the bounding 4--manifold. For
group--categories we evaluate the effect of these connected sums in terms of
structure constants of the category.
When the underlying group is cyclic the conclusions are very explicit, and
determine exactly when the field theory can
be normalized. For instance over an algebraically closed field there are
four categories with underlying group
$Z/2Z$,  distinguished by how  the non-unit simple object commutes with
itself. The possibilities are
multiplication by $\pm1$ or $\pm i$, $i$ a primitive fourth root of unity.
The $\pm1$ cases are symmetric, $\pm i$
braided--symmetric. The canonical and braided cases can be normalized; the
$-1$ case cannot.

\rk{Acknowledgement} The author is partially supported by the National 
Science Foundation.

\head Group--categories \endhead This section gives the formal
definition and examples, then proceeds to classification.
Classification is approached on three levels: modern homotopy theory
gives a quick general description.  Explicit CW models for classifying
spaces give associativity and commutativity isomorphisms satisfying
the standard axioms.  Finally chosing bases for morphism sets gives a
very explicit description in terms of sequences of units in the ring.
Much of this material is essentially already known, so proofs are
designed to clarify connections rather than nail down every detail.
For instance the iterated bar construction is explained in detail in
2.4 because the connection with categories comes from the details,
while the technically more powerful multi-simplicial construction
behind 2.3 is not discussed. We do give a lot of detail, though, since new
insights tend to be found in details.

\subhead 2.1 Definition\endsubhead A group--category is an additive category
over
 a commutative ring $R$, with a product (monoidal structure) that
 distributes over addition, and in addition:
\roster
\item it is additively semisimple in the sense that each object is a
finite sum of certain specified ``simple'' objects;
\item there are no nontrivial morphisms between distinct simple objects;
and \item the simple objects are invertible.  \endroster

An object is invertible if there is another object so that the product of
the two is isomorphic to the multiplicative unit.  This is a very
restrictive condition.  In particular it follows that the product of
any two simple objects is again simple, so isomorphism classes of simple
objects form a group.  This is called the ``underlying group'' of the
group category.  Condition (2) is usually automatic for simple objects
because the category is usually assumed to be abelian (have kernels
and cokernels, see \cite{McL}).  We avoid this assumption to enable use of
integers
and other non-fields as coefficient rings.  The extra generality is
useful in the abstract theory and really vital in some numerical
computations.

\subhead 2.2 Examples \endsubhead
The canonical examples are analogs of group rings.  Other examples
come from representations of groups and Lie algebras.
 \subsubhead 2.2.1
Canonical examples \endsubsubhead Suppose $G$ is a group and $R$ a
commutative ring.  Define $\Cal R[G]$ to be the category with objects
$G$--graded free $R$--modules of finite total dimension.  Morphisms are
$R$--homomorphisms that preserve the grading, with the usual
composition.  The product is the standard graded product: if $a$ and
$b$ are $G$--graded modules then $$(a\otimes b)_f =
\oplus_{\{g,h\: gh=f\}}a_g\otimes b_h.$$
Products of morphisms are defined similarly.  This
product is naturally associative with associating isomorphism the
``identity''
$$\oplus_{\{g,h,i\: (gh)i=f\}}(a_g\otimes b_h)\otimes c_i
\quad = \quad
\oplus_{\{g,h,i\: g(hi)=f\}}a_g(\otimes b_h\otimes c_i).$$
The simple objects are the ``delta functions'' that take all but one group
element to zero, and that one to a copy of $R$. In 2.3.3 and 2.4.1 we see
that general
group--categories are obtained  (up to equivalence) by modifying the
associativity and commutativity structures in this
standard example.

\subsubhead 2.2.2 Sub group--categories \endsubsubhead If $\Cal C$ is
an additive category with a product then the subcategory generated by
the invertible objects is a group--cateory.  The following examples  are of
this type.

\subsubhead 2.2.3 One--dimensional representations \endsubsubhead
If $R$ is a commutative ring and $G$ is a group then a representation of
$G$ over $R$ is a finitely generated free $R$--module on which $G$
acts.  Equivalently, these are $R[G]$ modules that are finitely
generated free as $R$--modules.  Tensor product over $R$ gives a
monoidal structure on the category of finite dimensional representations.  The
invertible elements in this category are the one--dimensional
representations.  Therefore the subcategory with objects
sums of 1--dimensional representations is a group category.  In fact it
is equivalent to the canonical group--category $\Cal
R[\hom(G,\units(R)]$.  We briefly describe the equivalences between
the two descriptions since they are models for several other constructions.

A homomorphism $\rho\:G\to \units(R)$ determines a 1--dimensional
representation $R^{\rho})$, where elements $g$ act by multiplication by
$\rho(g)$.

An object in the group--category is a free
$\hom(G,\units(R))$--graded $R$--module, so associates to each
homomorphism $\rho$ a free module $a_{\rho}$. Take such an object to
the representation
$\oplus_{\rho}(a_{\rho}\otimes_{R}R^{\rho})$.   This
clearly extends to morphisms.  The canonical identification
$R^{\rho}\otimes R^{\tau} = R^{\rho\tau}$ makes this a monoidal
functor from the group--category to representations.

To go the other way suppose $V$ is a representation.  Define a
$\hom(G,\units(R))$--graded $R$--module by associating to each
homomorphism $\rho$ the space $\hom(R^{\rho},$ $V)$.  To give an object
in the group--category these must be finitely generated free modules.
This process therefore defines a functor on the subcategory of
representations with this property, and this certainly contains sums
of 1--dimensional representations.  Note this functor may not be
monoidal on its entire domain: there may be indecomposable modules of
dimension greater than 1 whose product has 1--dimensional summands.
However it is monoidal on the subcategory of sums of 1--dimensional
representations.  It is also easy to see it gives an inverse
equivalence for the functor defined above.

\subsubhead 2.2.4 Quantum categories \endsubsubhead Let $G$ be a simple Lie
algebra, or more precisely an algebraic Chevalley group over $Z$, and
$p$ a prime larger than the Coxeter number of $G$. Some of the categories
are defined for non-prime $p$, but the prime
case is simpler and computations are currently limited to primes. Let
$X$ be the  weight lattice.  The ``quantum'' categories are obtained by:
consider
either mod $p$ representations (Gelfand and Kazhdan \cite{GK}), or
deform the universal enveloping algebra and then specialize the
deformation parameter to a $p^{\text{th}}$ root of unity [\Cite{L},
\Cite{Ka}].  Define $\Cal G_p$ to be the additive category generated by
highest weight representations whose weights lie in the standard
alcove of the positive Weyl chamber in $X$.  Define a product on $\Cal
G_p$ by: take the usual tensor product of representations and throw away
all indecomposable summands that are not of the specified type.  The
miracle is that this operation is associative, and gives a tortile or
symmetric monoidal category in the root of unity or mod $p$ cases
respectively.

Now let $R\subset X$ denote the root lattice of the algebra. The
quotient $X/R$ is a finite abelian group, and each highest weight
representation determines an element in $X/R$ (the equivalence
class of its weight).  Subgroups of $X/R$ correspond to
Lie groups with algebra $G$, and representations of the group are those with
weights in the given subgroup.  In particular the ``class 0''
representations, ones with weights in the root lattice, form a
monoidal subcategory.
\proclaim{ Conjecture} The category $\Cal G_p$
has a group subcategory with underlying group $X/R$, and $\Cal G_p$
decomposes as a tensor product of this subcategory and the class 0
representations $\Cal G_p^0$.  Further $\Cal G_p^0$ is ``simple'' in the sense
that it has no proper subcategories closed under products and
summands.  \endproclaim
This is true in the few dozen numerically computed examples, though
the tensor product in the root-of-unity cases might be slightly twisted.  In
these examples the objects in the group subcategory have weights lying
just below the upper wall of the alcove.  The values of these weights
are available through the ``Category Comparison'' software
in~\cite{QTP} (see the Category Guide).

\subhead 2.3 Homotopy classification of group--categories\endsubhead
Current homotopy-theory technology is used  to obtain the
classification in terms of spaces with two homotopy groups, or equivalently
group cohomology. The result is
essentially due to
[\Cite{FK}; \S7.5] where these are called ``$\Theta$--categories'':
 \proclaim{Proposition} Suppose
$R$ is a commutative ring and $G$ is a group.  Then \roster\item
monoidal group--categories over $R$ with underlying group $G$
correspond to $H^{3}(B_{G};\units(R))$; \item tortile (ie, balanced
braided--commutative monoidal) group--categories correspond to
$H^{4}(B^{2}_{G};\units(R))$; and \item symmetric monoidal
group--categories to $H^{d+2}(B^{d}_{G};\units(R))$, for $d>2$.
\endroster\endproclaim

It has been known for a long time that the group completion of the
nerve of a category with an associative monoidal structure is a loop
space.  It has been known almost as long that if the category is
symmetric then the group completion is an infinite loop space.
Recently this picture has been refined [\Cite{F}, \Cite{BFSV}] to include braided
categories:
the group completion of the nerve of a braided--commutative monoidal
category is a 2--fold loop space.  This can be applied to
group--categories to obtain:

\proclaim{2.3.1 Lemma} Suppose $R$ is a commutative ring and $G$ a group
(abelian in cases 2 and 3).
Then \roster \item equivalence classes of monoidal group--categories
over $R$ with group $G$ correspond to homotopy classes of simple
spaces with loop space\nl $B_{\units(R)}\times G$; \item
braided--commutative group categories correspond to spaces with second
loop space $B_{\units(R)}\times G$; and \item symmetric
group--categories correspond to spaces with $d$--fold loop\nl
$B_{\units(R)}\times G$, for $d>2$. \endroster \endproclaim
In practice this version seems to be more fundamental than the cohomology
description of the Proposition.
\demo{Proof} Consider the monoidal subcategory of simple objects and
isomorphisms in the category.  The nerve of a category is the simplicial
set with vertices the objects, and $n$--simplices for $n>0$ composable
sequences of morphisms of length $n$.  Condition 2.1(2) implies this
is a disjoint union of components, one for each isomorphism class.
Invertibility implies the components are all homotopy equivalent.
Endomorphisms of the unit object in a category over a ring $R$ are
assumed to be canonically isomorphic to $R$, so the isomorphisms of
each simple are given by $\units(R)$.  This identifies the nerve of
the whole category as $B_{\units(R)}\times G$.

The next step in applying the loop-space theory is group completion.
Ordinarily $\pi_{0}$ of a category nerve is a monoid, and group
completion converts this to a group.  Here $\pi_{0}$ is already the
group $G$ so the nerve is equivalent to its group completion.  Thus
application of [\Cite{F}, \Cite{BFSV}] shows that the nerve is a 1--, 2-- or
$d>2$--fold loop space when the category is monoidal,
braided--commutative, and symmetric respectively.

\noindent A few refinements are needed: \roster\item In the single loop
case the
delooping is $X$ with $\pi_{1}(X)=G$ and $\pi_{2}=\units(R)$.  In
general $\pi_{1}$ acts on higher homotopy groups.  Here the action is
trivial (the space is simple) because in the category $G$ acts
trivially on the coefficient ring.  \item It is not necessary to be
specific about which $d>2$ in the symmetric case because in this
particular setting a 3--fold delooping is automatically an infinite
delooping.  This follows from the cohomology description below.  \item
Generally the construction does not quite give a correspondence:
monoidal structures give deloopings of the group completion of the
nerve, while deloopings give monoidal structures on categories whose
nerve is already the group completion.  Here, however, the nerve is
group--completed to begin with, so the inverse construction does give
monoidal structures on categories equivalent to the original one.
\item Since the original group--category is additively semisimple,
monoidal structures on the simple objects extend linearly, and
uniquely up to equivalence, to products on the whole category that
distribute over sums.  This shows that classification of structures on
the subcategory of simples does classify the group--category.\qed
\endroster \enddemo

The final step in the classification is to relate this to group cohomology.

\proclaim{2.3.2 Lemma} Connected spaces with $\pi_{d}=G$,
$\pi_{d+1}=\units(R)$, and all other homotopy trivial (and simple if
$d=1$) are classified up to homotopy equivalence by elements of
$H^{d+2}(B^{d}_{G};\units(R))$.  \endproclaim

\demo{Proof} This is an almost trivial instance of Postnikov systems
[\Cite{Wh}; chapter {IX}].
Suppose $E$ is the space with only two non-vanishing homotopy groups.
There is a map $E\to B^{d}_{G}$ (obtained, for instance, by killing
$\pi_{d+1}$), and up to homotopy this gives a fibration
$$B^{d+1}_{\units(R)}\to E\to B^{d}_{G}.$$
The point of Postnikov systems is that this extends to the right:
there is a map $k\:B^{d}_{G}\to B^{d+2}_{\units(R)}$ well-defined up to
homotopy, so that
$$E@>>> B^{d}_{G}@>k>> B^{d+2}_{\units(R)}$$
is a fibration up to homotopy. This determines $E$, again up to homotopy.
Homotopy classes of such
maps $k$ are exactly
$H^{d+2}(B^{d}_{G};\units(R))$, so the spaces $E$ correspond to
cohomology classes.\qed\enddemo

Putting 2.3.1 and 2.3.2 together gives the classification theorem.

\subsubhead 2.3.3 Monoidal categories from spaces with two homotopy
groups\endsubsubhead
 In many ways the delooping of 2.3.1 is  more
fundamental than its $k$--invariant of 2.3.2. We finish this section by
showing how to recover the category from the
space. A  description directly in terms of the $k$--invariant is given in~2.4.
Suppose $E$ is a space with  $\pi_1(E)=G$, $\pi_2(E)=\units(R)$, and
$\pi_1$ acts trivially on $\pi_2$. This data
specifies (up to monoidal equivalence) a group--category over $R$ with
underlying group $G$. Here we show how to describe a category in the
equivalence class. In 2.3.4 this is extended to braided--monoidal and
symmetric categories.

Begin with the canonical category $\Cal R[G]$ of 2.2.1. $\Cal G$ has the
same underlying additive category over $R$
and the same product functor, but we change the associativity isomorphisms.
Specifically we find $\alpha(f,g,h)$ so
that the isomorphism $(a_f\otimes b_g)\otimes c_h\to a_f\otimes (b_g\otimes
c_h)$ obtained by multiplying the
standard isomorphism by $\alpha$ gives an associativity. The key property
is the pentagon axiom.

The definition of $\alpha$ depends on lots of choices. For each $g\in
\pi_1(E)$ choose a map $\hat g\:I/\partial
I\to E$ in the homotopy class. For each $g, h\in G$ choose a homotopy
$m_{g,h}\:\hat g\hat h\sim \widehat{gh}$. Here
$\hat g\hat h$ indicates composition of paths. The only restrictions are
that the identity element of the group
lifts to the constant path, and $m_{1,g}$ and $m_{g,1}$ are constant
homotopies.

Now define $\alpha(f,g,h)$ as follows: use these standard homotopies to
construct a homotopy $\widehat{fgh}\sim
\widehat{fg}\hat h\sim \hat f\hat g\hat h\sim \hat f\widehat{gh}\sim
\widehat{fgh}$. Since this is a homotopy of a
loop to itself the ends can be identified to give a map $I\times
S^1/(\partial I\times S^1)\to E$. Think of $I\times
S^1/(\{0\}\times S^1)$ as $D^2$, then this defines an element in
$\pi_2(E)=\units(R)$. Define $\alpha(a,b,c)$ to be
this element of $R$.

We explain why the pentagon axiom holds. A huge diagram goes with this
explanation, but the reader may find it
easier to reconstruct the diagram than to make sense of a printed version.
Thus we stick with words. It is sufficient
to verify the axiom for simple objects, and we write
$g$ for the
$G$--graded $R$--module that takes $g$ to $R$ and all other elements to 0.
The pentagon has various associations of a
4--fold product $efgh$ at the five corners, and connects them with
reassociation isomorphisms $\alpha$. The routine
for constructing the isomorphisms can be described as follows. Put the loop
$\widehat{efgh}$ at each corner, and put
the composite loop $\hat e\hat f\hat g\hat h$ in the center. Along each
radius from a corner to the center put the
concatenation of homotopies $m_{*,*}$ corresponding to the way of
associating the product at that corner. The
$\alpha$ for an edge comes from the homotopy of $\widehat{efgh}$ to itself
obtained by going from one
corner radially in to the center and then back out to  the other corner.
Going all the way around the pentagon
corresponds to going in and out five times. But going out and back in along
a single radius gives the composition of
a homotopy with its inverse, so cancels, up to homotopy. Therefore the
homotopy obtained from the full circuit is
homotopic to the constant homotopy of $\widehat{efgh}$ to itself. In
$\pi_2E=\units(R)$ this is the statement that
the product of the $\alpha$ terms associated to the edges is the identity,
so the diagram commutes.

Changing the choices gives an isomorphic category. Specifically suppose
$m'_{f,g}$ are different homotopies between
compositions. They differ from the original $m$ by elements of $\pi_2(E)$,
so by units $\tau_{f,g}\in R$. Regard
this as defining a natural isomorphism from the product functor to itself:
$f\otimes g\to f\otimes g$ by
multiplication by $\tau_{f,g}$. Then the identity functor $\Cal R[G]\to\Cal
R[G]$ together with this transformation
is a monoidal isomorphism, ie, associativity defined using
$m$ in the domain commutes with associativity using $m'$ in the range. We
revisit this construction in the context
of group cohomology in 2.4.2, and make it more explicit using special
choices in 2.5.

\subsubhead 2.3.4 Braided group--categories from spaces with two homotopy
groups\endsubsubhead
Suppose $E$ has $\pi_2E=G$ and $\pi_3E=\units(R)$. According to 2.3.1 this
corresponds to an equivalence class of
braided--commutative group--categories with underlying group $G$. Here we
show how to extract one such category from
this data, extending the monoidal case of 2.3.3.

The loop space $\Omega E$ has $\pi_1E=G$ and $\pi_2E=\units(R)$, so
specifies an associativity
structure for the standard product on $\Cal R[G]$. Let $\{\hat g\}\:I\to
\Omega E$ and $m_{f,g}\:I^2\to \Omega E$ be
the choices used in 2.3.3 to make this explicit. Let $\widetilde g\:I^2\to
E$ and $\widetilde m_{f,g}\:I^3\to E$ denote the
adjoints. An element $\sigma(f,g)\in\units(R)$ is obtained as follows:
define a homotopy $\widetilde{fg}\sim \widetilde
f\widetilde g\sim \widetilde g\widetilde f\sim \widetilde{gf}$ by the
reverse of $\widetilde m_{f,g}$, the clockwise standard commuting
homotopy in $\pi_2$, and $\widetilde m_{g,f}$. Since $G$ is abelian
$gf=fg$, and this is a self-homotopy. Glueing the
ends gives a map on $I^2\times S^1/(\partial I^2\times S^1)$. Regard this
as a neighborhood of $S^1\subset D^3$, and
extend the map to $D^3$ by taking the complement to the basepoint. This
gives an element of $\pi_3(E)=\units(R)$.
Define this to be $\sigma(f,g)$. Define a commutativity natural
transformation $f\otimes g\to g\otimes f$ by
multiplying the natural identification by $\sigma(f,g)$.

We explain why this and the associativity from 2.3.3 satisfy the hexagon
axiom. Again we omit the huge diagram. The
hexagon has various associations of permutations of $fgh$ at the corners,
and reassociating and commuting
isomorphisms alternate going around the edges. Imagine a triangle inside
the hexagon, with two hexagon corners
joined to each triangle corner. Put
$\widetilde{fgh}$ at each hexagon corner, and the three permutations of
$\widetilde f\widetilde g\widetilde h$ on the triangle
corners. On the edges joining the triangle to the hexagon put compositions
of homotopies $\widetilde m$ corresponding to
different associations of the terms. On the edges of the triangle put
clockwise commuting homotopies in $\pi_2$. The
homotopies used to define associating or commuting units on the hexagon
edges are obtained by going in to the
triangle and either directly back out (for associations) or along a
triangle edge and back out (for commutes). Going
around the whole hexagon composes all these. The trips from the triangle
out and back cancel, to give a homotopy of
the big composition to the composition of the triangle edges. This
composition is trivial (it gives the analog
of the hexagon axiom for $\pi_2(E)$). Thus the composition of homotopies
corresponding to the full circuit of the
hexagon gives the trivial element in $\pi_3(E)$, and the diagram itself
commutes.

Finally suppose $E$ has $\pi_d(E)=G$ and $\pi_{d+1}(E)=\units(R)$ for some
$d>2$. Then the same arguments as above
apply except the representatives $\widetilde g$ are now defined on $D^d$,
and there is only a single standard commuting
homotopy, up to homotopy. This implies $\sigma(g,h)\sigma(h,g)=1$, so the
group--category is symmetric.

\subhead 2.4 Models for classifying spaces \endsubhead
Here we use explicit CW models for classifying spaces $B^{n}_{G}$ to
connect group cohomology to descriptions of categories
using functorial isomorphisms.  This is done in detail for $n=1,2$,
and outlined for $n=3$.  The basis for the connection is a comparison
between general group--categories and the standard example~2.2.1.

\proclaim{2.4.1 Lemma} Suppose $\Cal C$ is a group category over $R$
with underlying group $G$.  Then there is an equivalence of categories
$\Cal C\to \Cal R[G]$ and a natural transformation between the given
product in $\Cal C$ and the standard product in $\Cal R[G]$.\endproclaim

Note that this functor usually not monoidal since it usually will not
commute with associativity morphisms.  If a ``lax'' description of
associativity is used  then it can be transferred through such a
categorical equivalence.  The classification of group--categories then
corresponds to classification of different associativity and
commutativity structures for the standard product on $\Cal R[G]$.

\demo{Proof} By hypothesis $G$ is identified with the set of
equivalence classes of simple objects in $\Cal C$, so we can choose a
simple object $s_{g}$ in each equivalence class $g$.  Further we can
choose isomorphisms $m_{g,h}\:s_{gh}\to s_{g}\sq s_{h}$.

Now define the functor $\hom_{s}\:\Cal C\to \Cal R[G]$ by: an object
$X$ goes to the function that takes $g\in G$ to $\hom_{\Cal
C}(s_{g},X)$.  Comparison of products in the two categories involves
the diagram $$\CD \Cal C\times \Cal C @>{\hom_{s}\times\hom_{s}}>>\Cal
R[G]\times \Cal R[G]\\
@VV{\sq}V @VV{\otimes}V\\
\Cal C @>\hom_{s}>> \Cal R[G]\endCD$$
A natural transformation $\otimes (\hom_{s}\times\hom_{s})\to
\hom_{s}\sq $ consists of: for $X, Y$ in $\Cal C$ and $g\in G$ a natural
homomorphism $\oplus_{h}\hom(s_{h},X)\otimes\hom(s_{h^{-1}g},Y)\to
\hom(s_{g},X\sq Y)$.  Define this by taking
$(a,b)\in \hom(s_{h},X)\otimes\hom(s_{h^{-1}g},Y)$ to $(a\sq
b)m_{h,h^{-1}g}$.  It is simple to check this has the required
naturality properties.  Note the lack of any coherence among the
isomorphisms $m_{g,h}$ prevents any conclusions about associativity. \qed\enddemo

Associativity structures for a product on a category are defined using
natural isomorphisms satisfying the ``pentagon axiom'' \cite{McL}.
These can be connected directly to group cohomology via the cellular
chains of a particular model for the classifying space.

\proclaim{2.4.2 Lemma} Suppose $G$ is a group and $R$ a commutative ring.
\roster \item Cellular 3--cocycles for the bar construction $B_{G}$ are
natural associativity isomorphisms for the product on $\Cal R[G]$,
and coboundaries of 2--cochains correspond to compositions with natural
endomorphisms.  \item If $G$ is abelian, cellular 4--cocycles for the
iterated bar
construction $B^{2}_{G}$ give braided--commutative monoidal structures for
the
product on $\Cal R[G]$, and coboundaries of 3--cochains correspond to
natural endomorphisms.
\item If $G$ is abelian, cellular 5--cocycles on $B^{3}_{G}$ give symmetric
monoidal structures.
\endroster\endproclaim

Lemmas 2.4.1 and 2.4.2 together give the equivalences between
group--categories and cohomology, except for ``balance'' in the braided
case.  This is addressed in 2.4.3.  The analysis in the symmetric case
is only sketched.

\demo{Proof} Suppose $G$ is a discrete group.  The ``bar construction''
gives the following model for the classifying space $B_{G}$: $n$--cells
are indexed by $n$--tuples $(g_{1},\dots,g_{n})$ of elements in the
group, so we denote the set of $n$--tuples by $B^{(n)}_{G}$.  Note that
there is a single 0--cell, the 0--tuple $(\;)$.  There are $n+1$
boundary functions from $n$--tuples to $(n-1)$--tuples: $\partial_{0}$
omits the first element; $\partial_{n}$ omits the last; and for
$0<i<n$, $\partial_{i}$ multiplies the $i$ and $i+1$ entries:
$\partial_{i}(g_{1},\dots,g_{n}) =
(g_{1},\dots,g_{i}g_{i+1},\dots,g_{n})$.

 We get a space by geometrically
realizing these formal cells:
$$B_{G}=\bigl(\cup_{n}B^{(n)}_{G}\times \Delta^{n}\bigr)/\simeq.$$
Here $\Delta^{n}$ is the standard $n$--dimensional simplex, and
$\simeq$ is the equivalence relation that for each $n$--tuple $\tau$
identifies $\tau\times\partial_{i}\Delta^{n}$ with
$\partial_{i}\tau\times\Delta^{n}$.

The cellular chains of this CW structure gives a model for the chain
complex of the space.  Specifically, $C^{c}_{n}(B_{G})$ is the free
abelian group generated by the formal $n$--cells $B^{(n)}_{G}$, and the
boundary homomorphism $\partial\:C^{c}_{n}(B_{G})\to
C^{c}_{n-1}(B_{G})$ takes an $n$--tuple $\tau$ to the class
representing the boundary $\partial \tau$.  The boundary of the standard
$n$--simplex $\Sigma^{n}$ is the union of the faces
$\partial_{i}\Sigma^{n}$, but the ones with odd $i$ have the wrong
orientation.  Using the equivalence relation in $B_{G}$ therefore gives
$\partial\tau= \Sigma_{i=0}^{n}(-1)^{i}\partial_{i}\tau$.

Now suppose $H$ is an abelian group.  The model for chains of $B_{G}$
gives a description for the cohomology $H^{3}(G;H)$.  A 3--cocycle is a
function $\alpha\:B^{3}_{G}\to H$ with composition $\alpha\partial$ is
trivial.  $\partial(a,b,c,d) =
(b,c,d)-(ab,c,d)+(a,bc,d)-(a,b,cd)+(a,b,c)$, so the cocycle condition
is
$$\alpha(b,c,d)+\alpha(a,bc,d)+\alpha(a,b,c)=\alpha(ab,c,d)+\alpha(a,b,cd).$$
In the application the coefficient group is $\units(R)$, with
multiplication as group structure.  Rewriting the cocycle condition
multiplicatively gives exactly the pentagon axiom for associativity,
so this gives a monoidal category.

Now we consider uniqueness. A 2--cochain is a function on the 2--cells, so
$\mu(a,b)$ defined for all $a,b\in G$.  The coboundary of this is the
3--cochain obtained by composing with the total boundary homomorphism.
Written multiplicatively (in $\units(R)$) this is
$$(\delta \mu)(a,b,c)=\mu(b,c)\mu(ab,c)^{-1}\mu(a,bc)\mu(a,b)^{-1}.$$
Thus a 3--cocycle $\alpha'$ differs from $\alpha$ by a coboundary if
$$\alpha'(a,b,c)=\mu(a,b)^{-1}\mu(ab,c)^{-1}\alpha(a,b,c)\mu(b,c)\mu(a,bc).$$
Interpreting this as commutativity in the diagram
$$\CD (ab)c @>\alpha(a,b,c)>>(a(bc)\\
@VV{\mu(a,b)}V @VV{\mu(b,c)}V\\
(ab)c @.  (a(bc)\\
@VV{\mu(ab,c)}V @VV{\mu(a,bc)}V\\
(ab)c @>\alpha'(a,b,c)>>(a(bc)\endCD$$
shows we can think of $\mu$ as a natural transformation from the
standard product to itself, and then $\alpha'$ is obtained from
$\alpha$ by composition with this transformation.  This gives an
isomorphism between categories where the associativity cocycles differ by
a coboundary.

The braided case uses  the
iterated bar construction. If $G$ is abelian then $B_{G}$ is again a group,
this time simplicial or topological rather than discrete.  The same
construction gives a simplicial (or $\Delta$) space $B(B_{G})$ whose
realization is $B^{2}_{G}$.  The first step in describing this is a
description of the multiplication on $B_{G}$.

Cells in the product
$B_{G}\times B_{G}$ are modeled on products
$\Delta^{i}\times\Delta^{j}$.  The map $B_{G}\times B_{G}\to B_{G}$ is
defined by subdividing these products into simplices, and describing
where in $B_{G}$ to send these simplices.
The standard subdivision of a product of simplices is obtained as follows: the
vertices of $\Delta^{i}$ are numbered $0, 1,\dots, i$.  Suppose
$((r_{0},s_{0}),\dots, (r_{i+j},s_{i+j}))$ is a sequence of pairs
of these, ie,  vertices of $\Delta^{i}\times\Delta^{j}$, then the
function of vertices $k\mapsto (r_{k},s_{k})$ extends to a linear map to the
convex hull $\Delta^{i+j}\to \Delta^{i}\times\Delta^{j}$.  Restrict
the sequences to ones for which one coordinate of $(r_{k+1},s_{k+1})$ is
the same as in $(r_{k},s_{k})$, and the other coordinate increases by
exactly one.  Then this gives a collection of embeddings with disjoint
interiors, whose union is the whole product.

We relate this subdivision to the indexing of simplices by sequences
in $G$.  Think of a sequence $(a_{1},\dots,a_{i})$ as labeling edges
in $\Delta^{i}$, specifically think of $a_{k}$ as labeling the edge
from vertex $k-1$ to $k$.  Then we label sub-simplices of a product
$(a_{*})\times (b_{*})$ by: if $r_{k}=r_{k-1}+1$ then label the edge
from $k-1$ to $k$ with $a_{r_{k}}$, otherwise label it with
$b_{s_{k}}$.  This identifies the sub-simplices as corresponding to
$i,j$--shuffles: orderings of the union $(a_{*})\cup(b_{*})$ which
restrict to the given orderings of $a_{*}$ and $b_{*}$.  Thus we can write
$$\Delta^{i}\times\Delta^{j} = \cup_{s}s(\Delta^{i+j})$$
where the union is over $i,j$ shuffles $s$.  For future reference we
mention that the orientations don't all agree: the orientation on
$s(\Delta^{i+j})$ is $(-1)^{s}$ times the orientation on the product,
where $(-1)^{s}$ indicates the parity of $s$ as a permutation.

Now the product on $B_{G}$ is defined by: if $s$ is a shuffle the
sub-simplex $s(\Delta^{i+j})\times(a_{*})\times(b_{*})$ goes to
$\Delta^{i+j}\times s(a_{*},b_{*})$.  It is a standard fact that this is
well-defined on intersections of sub-simplices.

As before the $n$--simplices of $B(B_{G}$ are indexed by points in the
$n$--fold product $\times^{n}B_{G}$.  The realization is again
$$B^{2}_{G} = \bigl(\cup_{n}\Delta^{n}\times (\times^{n}B_{G})\bigr)/\simeq.$$
The equivalence relation identifies points in the boundary of
$\Delta^{n}$ with points in lower-dimensional pieces.  Specifically we
identify
$\partial_{k}\Delta^{n}\times (\times^{n}B_{G})$ with its image in
$\Delta^{n-1}(\times^{n-1}B_{G})$, via the map which is the
``identity'' on the simplices, and on the $B_{G}$ part multiplies
the $k-1$ and $k$ entries if $0<k<n$, omits the first if $0=k$, and omits
the last if $k=n$.

This definition gives a cell complex model for $B^{2}_{G}$.
Unraveling, we find the cells are of the form
$$\Delta^{n}\times\bigl(\Delta^{i_{1}}\times\cdots\times\Delta^{i_{n}}\bigr)
\times \bigl((a^{1}_{*})\times\cdots \times (a^{n}_{*})\bigr),$$
where $(a^{k}_{*})$ is a sequence of length $i_{k}$.

The cell structure on the space gives standard models for the chain
and cochain complexes.  The first comment about the chain complex is
that the cells that involve the 0--cell of $B_{G}$ form a contractible
subcomplex.  The union is not a topological subcomplex because these
cells  have faces that are not of this type.  However if a face
does not involve a 0--cell then there is an adjacent face with the same
image but opposite sign, so they algebraically cancel in the chain
complex.  Dividing out this subcomplex leaves  ``non-trivial'' cells,
corresponding to non-empty sequences $(a_{*})$.

We use this to describe the cohomology group $H^{4}$.  Eventually the
coefficients will be $\units(R)$, but to keep the notation standard we
start with a group $J$ with group operation written as addition.
Nontrivial 4--cells are in two families: $$\align
\Delta^{1}\times(\Delta^{3})&\text{ indexed by }(a,b,c) \text{,
and}\\\Delta^{2}\times(\Delta^{1}\times\Delta^{1}) &\text{ indexed by
}((a),(b))\endalign$$
Denote the cochain $C_{4}\to J$ by $\alpha(a,b,c)$ on the first
family, and $\sigma(a,b)$ on the second.

The cocycle condition on $(\alpha, \sigma)$ comes from boundaries of
5--cells.  Nontrivial 5--cells are in families:
$$\align
\Delta^{1}\times(\Delta^{4})&\text{ indexed by }(a,b,c,d)\\
\Delta^{2}\times(\Delta^{1}\times\Delta^{2}) &\text{ indexed by
}((a),(b,c)) \text{,
and}\\
\Delta^{2}\times(\Delta^{2}\times\Delta^{1}) &\text{ indexed by
}((a,b),(c))\endalign$$
In the first family the boundary of the $\Delta^{1}$ factor is trivial, so the
boundary is the boundary of $(a,b,c,d)$ as a 4--cell of $B_{G}$.  As
before this gives the pentagon axiom for $\alpha$.
Now consider $((a),(b,c))$ in the second family.  Boundaries of
products are given by $\partial(x\times y)=\partial(x)\times y
+(-1)^{\text{dim\,}(y)}x\times\partial y$.  In
$\Delta^{2}\times(\Delta^{1}\times\Delta^{2})$ the boundary on the
middle piece vanishes so the total boundary is
$\partial\times\id\times\id -\id\times\id\times\partial$.  In the
first factor the boundary is $\partial_{0}-\partial_{1}+\partial_{2}$.
The first and last use projection of $B_{G}\times B_{G}$ to
one factor, so map to cells of dimension less than 3 and are trivial
algebraically.  $\partial_{1}$ uses multiplication in
$B_{G}$ so is given by $(1,2)$ shuffles.  This contribution to the
boundary is thus $-\bigl((a,b,c)-(b,a,c)+(b,c,a)\bigr)$.  The boundary
in the last coordinate applies the $B_{G}$ boundary to $(b,c)$.  This
contribution is $-\bigl(((a),(c))-((a),(bc))+((a),(b))\bigr)$.  Applying the
cochain and setting it to zero gives
$$\alpha(a,b,c)-\alpha(b,a,c)+\alpha(b,c,a)
+\sigma(a,c)-\sigma(a,bc)+\sigma(a,b) =0.$$
This is exactly the hexagon axiom for $\alpha$ and $\sigma^{-1}$,
written additively.  Boundaries of cells in the third family give the
hexagon axiom for $\alpha$ and $\sigma$.

The conclusion is that 4--dimensional cellular cochains in $B^{2}_{G}$
correspond exactly to associativity and commutativity isomorphisms
$(\alpha,\sigma)$ satisfying the pentagon and hexagon axioms, for
the standard product on the category $\Cal R[G]$.

The final step in the proof of Lemma 2.4.2 is seeing that coboundaries
correspond to endomorphisms, or more precisely natural transformations
of the standard product to itself.

The only nontrivial 3--cells in $B^{2}_{G}$ are of the form
$\Delta^{1}\times\Delta^{2}\times(a,b)$.  3--cochains therefore
correspond to functions $\mu(a,b)$.  Boundaries of 4--cells are given by:
in the $\Delta^{1}\times \Delta^{3}$ case, the negative of the $B_{G}$
 boundary (the negative comes from the preceeding $\Delta^{1}$ factor).
 In the $\Delta^{2}\times(\Delta^{1}\times\Delta^{1})$ case all terms
 vanish except the $-\partial_{1}$ term in the first factor, which gives
 shuffles $-\bigl((a,b)-(b,a)\bigr)$.  Therefore changing a 4--cocycle
$(\alpha,\sigma)$ by the coboundary of $\mu$ changes $\alpha(a,b,c)$
just as in the monoidal category case, and changes $\sigma$ by
conjugation by $\mu$.

Finally we come to the symmetric monoidal case, using $B^{3}_{G}$.
This is a further bar construction obtained as
$$B^{3}_{G}=\bigl(\cup_{n}(\Delta^{n}\times(\times^{n}B^{2}_{G}))\bigr)/\simeq$$
where the identifications in $\simeq$ involve a product structure on
$B^{2}_{G}$.  We indicate the source of the new information (symmetry
of $\sigma$) without going into details.

We are concerned with $H^{5}$, so functions on the 5--cells.  Again we
can divide out the ``trivial'' ones involving 0--cells of $B_{G}$ at
the lowest level.  The only nontrivial cells are products of
$\Delta^{1}$ and 4--cells of $B^{2}_{G}$, so these use the same data
$(\alpha, \sigma)$ as 4--cochains on $B^{2}_{G}$.  Boundaries of
6--cells of the form $\Delta^{1}$ times a 5--cell of $B^{2}_{G}$ involve
only the second factor, so give the same relations as in $B^{2}$
(namely, the pentagon and hexagon axioms).  The only other source of
relations are nontrivial 6--cells of the form
$\Delta^{2}\times((\text{2--cell})\times(\text{2--cell}))$, where each of
these 2--cells (in $B^{2}$) is of the form
$\Delta^{1}\times\Delta^{1}_{(a)}$.  The only nonzero term in the
boundary of such a 6--cell comes from $\partial_{1}$ in $\Delta^{2}$,
which goes to $\Delta^{1}$ times the product of the two 2--cells in
$B^{2}$.  We won't describe this in detail, but multiplying two cells
of the form $\Delta^{1}\times\Delta^{1}$ involves multiplying the
first two $\Delta^{1}$ factors to get a square, then subdividing this
into two $\Delta^{2}$.  These two sub-simplices have opposite
orientation, so the product is a difference of cells of $B^{2}$ of the
form $\Delta^{2}\times(\Delta^{1}\times \Delta^{1})$.  Vanishing of
the cocycle on this therefore is a relation of the form  (written additively)
$\sigma(*)-\sigma(**)=0$.  So much follows from generalities.
We don't do it here, but explicit description of the product structure
shows the indices on the two final $\Delta^{1}$ factors is
interchanged, so we get exactly the symmetry relation
$$\sigma(a,b) = \sigma(b,a).\eqno{\qed}$$
\enddemo

\subhead 2.5 Numerical presentations\endsubhead
Here we get explicit  ``numerical
presentations''  of group--categories in the sense of \cite{BQ}.  This
amounts to direct
computation of group cohomology, and we interpret some of the formulae in
terms of cohomology operations.  We consider the symmetric and
braided--commutative cases in detail,
and only remark on the general monoidal case.

\proclaim{2.5.1 Proposition} Suppose $G$ is an abelian group with generators
$g_i$
of order $n_i$, and $R$ is a commutative ring.
\roster\item Braided--commutative group--categories over $R$ with underlying
group
$G$ correspond to

{\rm i)}\qua $\sigma_*$ with $\sigma_i^{2n_i}=1$, and $\sigma_i^{n_i}=1$
if $n_i$ is odd; and

{\rm ii)}\qua $\sigma_{i,j}$ for $i>j$, with
$\sigma_{i,j}^{n_i}=\sigma_{i,j}^{n_j}=1$.
\item These categories are all tortile, and any tortile structure is
obtained by scaling a standard one by a homomorphism
from $G$ to the units of $R$.
\item The symmetric monoidal categories correspond to
$\sigma_i^2=\sigma_{i,j}=1$.
\endroster
\endproclaim

Given a group--category we extract the invariants as follows:
Choose a simple object $\hat g_i$ in the equivalence class $g_i$.  The
commuting isomorphism $\sigma_{\hat g_i,\hat g_i}$ is an endomorphism of
the object $\hat g_i\sq \hat g_i$, so is multiplication by an element of
$R$.  Define this to be $\sigma_i$.  If $i>j$ the double commuting
isomorphism $\hat g_i\sq\hat g_j\to \hat g_j\sq\hat g_i\to \hat g_i\sq\hat
g_j$ is also an endomorphism, so is multiplication by an element of $R$.
Define this to be $\sigma_{i,j}$.

Conversely given invariants we define a group--category by defining
associating and commuting isomorphisms for the standard product on the
standard group--category $\Cal R[G]$.  The content of 2.5.1 is then that
there is a braided--monoidal equivalence from a general group--category to
the standard one with the same invariants.

\subsubhead 2.5.2 The inverse construction\endsubsubhead Suppose data as in
2.5.1 are given. If $a$ is an element
of $G$ we let $a_i$ denote the exponent of $g_i$ in $a$, so we have
$a=\Pi_ig_i^{a_i}$.  Then define: $$\align \alpha\:(ab)c\to a(bc)&
\text{ is multiplication by } \Pi_i\cases 1&\text{ if }b_i+c_i<n_i\\
\sigma_i^{n_ia_i}&\text{ if }b_i+c_i\geq n_i\endcases \\
\sigma\:ab\to ba&\text{ is multiplication by } \Pi_{i\leq
j}\sigma_{i,j}^{a_ib_j}
\endalign$$
In the second expression $\sigma_{i,j}$ means $\sigma_i$ if $i=j$.  Recall
the exponent of $\sigma_i$ is at most $2n_i$, so the terms
$\sigma_i^{n_ia_i}$ are 1 if $n_i$ is odd, and depend at most on the parity
of $a_i$ in general.

\subsubhead 2.5.3 Example \endsubsubhead Suppose $G=Z/2Z$ and $i\in R$ is a
primitive $4^{th}$ root of unity. Then
$\sigma=\pm 1$ and $\sigma=\pm i$  give four  group--categories that are not
braided commutative equivalent. The $\pm1$
cases are symmetric, and monoidally equivalent (ignoring commutativity). The
$\pm i$ cases are genuinely braided. In these the associativity $(gg)g\to
g(gg)$ is multiplication by $-1$, so they are
monoidally equivalent to each other but distinct from the standard category.

\subsubhead 2.5.4 A relation to cohomology \endsubsubhead Since
group--categories correspond to cohomology classes, Proposition 2.5.1
amounts to an explicit calculation of cohomology.  We discuss only a
piece of this: the associativity structure is the image of the braided
structure
under the suspension $$\Sigma\:H^4(B^2_G;\units(R))\to H^3(B_G;\units(R)).$$
Elements of $H^3(B_G;\units(R))$ can be obtained as follows:
\roster\item take homomorphisms $G\to J\to Z/2\to \units(R)$, with $J$ cyclic;
\item the identity homomorphisms defines a class $\iota\in H^1(B_J;J)$;
\item the Bockstein is an operation $\beta\:H^1(B_J;J)\to H^2(B_J;J)$;
\item applying the Bockstein to $\iota$ and then cup product with $\iota$
gives $\iota\cup\beta(\iota)\in H^3(B_J;J)$;
\item applying $B_G\to B_J$ in the space argument, and $J\to \units(R)$ in
the coefficients gives an element in $H^3(B_G;\units(R))$.
\endroster
Working out  the Bockstein and cup product  on the chain level gives
exactly the formulas in the description of $\alpha$ above when
$\sigma_i^{n_i}\neq 1$.

\subsubhead{2.5.5 Representatives and products}\endsubsubhead To begin
the construction we need: \roster\item a standard representative for
each isomorphism class of simple object; and \item an algorithm for
finding a parameterization of an arbitrary iterated product by the
standard representative.  \endroster Here we will use the solution to
the word problem in the abelian group $G$.  The analysis of other
categories uses the same approach, as far as it can be taken.
Descriptions of representations of $sl(2)$, cf  \cite{CFS}, and
other small algebras \cite{Ku} depend on the description of specific
representatives for simples using projections on iterated products of
``fundamental'' representations.  When special information of this
type is not available numerical presentations can be obtained by
numerically describing representatives and then parameterizing
iterated products by direct computation~[\Cite{B}, \Cite{BQ}].

Choose representatives as follows: choose simple objects $\hat g_i$ in the
equivalence class of the generator $g_i\in G$, for each $i$.  A general
element $a\in G$ has a unique representation of the form
$a=g_1^{r_1}g_2^{r_2}\cdots g_k^{r_k}$, where $0\leq r_i <n_i$.  We want to
get an object in the category by subsituting the simple object $\hat g_i$
for the group element $g_i$, but for this to be well-defined we must
specify a way to associate the product.  Associate as follows: each
$g_i^{r_i}$ is nested left (ie,  $g^4=((g g)g)g$), and then the product of
these pieces is also nested left.  Now subsituting standard representatives
for generators  gives a standard simple object in each equivalence class.

Next fix for each $i$ an isomorphism $\Lambda_i\:1\to g_i^{n_i}$.  Suppose
$W$ is a word with associations, in the generators $g_i$.  $W$ specifies an
iterated product, and we want an algorithm describing a morphism from the
standard representative for this simple object into the product of the word
$W$.  Proceed as follows:
\roster\item If there is a pair $g_jg_1$ with $j>1$ in the word (ignoring
associations), then associate to pair them, and
apply $\sigma_{g_1,g_j}^{-1}$ to interchange them. The result is a simpler
word $W'$ with a morphism (of products)
$W'\to W$ formed by composing associations and $\sigma_{g_i,g_j}$;
\item when (1) is no longer possible, then all $g_1$ occur first.  Repeat
to move all $g_2$ just after the $g_1$, etc.  Then associate to the left
to obtain $g_1^{r_1}g_2^{r_2}\cdots g_n^{r_n}$.
\item after (2) is done, if any $r_i$ is too large, compose with
$\Lambda_i\sq \id\:g_i^{r_i-n_i}\to g_i^{r_i}$.
\endroster
When this process terminates the result is a morphism from a standard
representative to the product of the word $W$.

\proclaim{ Lemma} The morphisms resulting from this algorithm are well
defined.
\endproclaim
The point is that there are choices, but the final result is
independent of these choices.  Suppose  that we have two sequences of
operations as described in the algorithm.  The coherence theorem for
associations shows the outcome does not depend on the order of
associations, so problems can come only from the $\Lambda$ in (3) and the
commuting isomorphisms in (1) and (2).  There is no choice about which
operations are needed, but some choice in the order.  If there is a choice
then the operations do not overlap, in the sense that each is of the form
$\id\sq\sigma\sq\id$, and the nontrival part of one operation takes place
in an identity factor of the other.  Thus the operations commute, and the
result is well-defined.

\subsubhead{2.5.6 The functor }\endsubsubhead We use the choices of 2.5.5
to define a functor $\Cal F\:\Cal C\to \Cal R[G]$, and a natural
transformation between the two products.

Suppose $a$ is an object of $\Cal C$.  $\Cal F(a)$ is supposed to be a
function from $G$ to $R$--modules.  Define $\Cal F(a)(g)=\hom_{\Cal C}(\hat
g,a)$, where $\hat g$ denotes the standard simple object in the equivalence
class $g$.  The natural transformation from the product in $\Cal R[G]$ to
the one in $\Cal C$ is given by natural homomorphisms
$$\oplus_{\{r,s\mid rs=g\}}\hom(\hat r,a)\otimes\hom(\hat s,b)@>>>
\hom(\hat g,a\sq b).$$
These are defined by $h_1\otimes h_2\mapsto (h_1\sq h_2)m$, where
$m\:\widehat{rs}\to \hat r\sq \hat s$ is the standard parameterization,
ie,  the morphism from the standard representative of the product to the
product of representatives.

The proposition is proved by showing this functor and transformation
commute with commutativity and
associativity isomorphisms when the twisted structure 2.5.2 is used in $\Cal
R[G]$.  we begin with very special cases.  Consider the commutativity
$\sigma_{\hat g_i,\hat g_j}\:\hat g_i\sq\hat g_j \to g_j\sq\hat g_i$.  If
$i<j$ then the right term is already canonical and the algorithm gives
$\sigma_{\hat g_j,\hat g_i}^{-1}$ as parameterization of the left.  The
diagram
$$\CD \hat g_i\sq\hat g_j @>\sigma_{\hat g_i,\hat g_j}>> g_j\sq\hat g_i\\
@AA{\id}A @AA{\sigma_{\hat g_j,\hat g_i}^{-1}}A\\
\widehat{g_ig_j}@>>>\widehat{g_ig_j}\endCD$$
commutes if we put $\sigma_{i,j}=\sigma_{\hat g_j,\hat g_i}\sigma_{\hat
g_i,\hat g_j}$ across the bottom.

If $i>j$ in the same situation then the left term is canonical and we get
the diagram
$$\CD \hat g_i\sq\hat g_j @>\sigma_{\hat g_i,\hat g_j}>> g_j\sq\hat g_i\\
@AA{\sigma_{\hat g_i,\hat g_j}^{-1}}A @AA{\id}A\\
\widehat{g_ig_j}@>>>\widehat{g_ig_j}\endCD$$
which commutes with the identity across the bottom.  The commutativity
required in the model is therefore multiplication by
 $$\cases \sigma_{i,j}&\text{ if } i<j\\
 \sigma_i&\text{ if } i=j\\
 1&\text{ if } i>j\endcases$$
 which is the factor specified in~2.5.2.

 Associativity terms come from different ways of reducing excessively large
 powers.  Fix a particular generator $g_i$, drop $i$ from the notation, and
 consider the association $(\hat g^r \hat g^s)\hat g^t\to \hat g^r( \hat
 g^s\hat g^t)$.  If $s+t<n$ then the parameterization algorithm gives the
 same thing on the two sides, and the associativity is the identity.  If
 $s+t\geq n$ the reductions using $\Lambda$ are different:
 $$\CD (\hat g^n)\sq \hat g^{r+s+t-n}@>\alpha >>\hat g^r\sq (\hat g^n)\sq \hat
 g^{s+t-n}\\
 @AA{\Lambda\sq \id}A @AA{\id\sq\Lambda\sq\id}A\\
 (1)\sq \hat g^{r+s+t-n}@>>>\hat g^r\sq (1)\sq\hat
 g^{s+t-n}\endCD$$
 Putting multiplication by
 $\sigma_{i}^{n_{i}r}$ on the bottom makes the diagram commute.  This
corresponds to commuting
  $g^{r}$ past $g^{n}$ one $g$ factor at a time.  The point is that this is
different
  from commuting the full products, which wouldn't contribute anything
  since $g^{n}=1$.

We now claim the hexagon axiom and these special cases imply the
general case, ie,  the associativity and commutativity isomorphisms in
$\Cal C$ commute with the natural transformation between products and
the twisted associativity and commutativity morphisms~2.5.2.  The new
feature in the general case is that different associations change the
way a product is reduced to standard form.  Specifically,
$g_{3}(g_{2}g_{1})$ follows the standard algorithm in first commuting
the $g_{1}$ all the way to the left, while in $(g_{3}g_{2})g_{1}$ the
$g_{3}g_{2}$ are commuted first.  However the fact that both orders
give the same final morphism is exactly the standard crossing identity
for braided--commutative categories.  Independence of association in
arbitrary products follows from this by induction on the number of
out-of-order commutes.  Once one can choose associations arbitrarily
it is straightforward to check the general associativity and
commutativity formulae by choosing special association patterns.

The general associative (ie,  non-braided) case of classification is
not considered in this section, but at this point we can indicate what
is involved when the underlying group is abelian.  As above choose
simple objects representing generators, and reduction isomorphisms
$\Lambda_{i}\:1\to \hat g_{i}^{n_{i}}$.  Since the underlying group is
abelian there are isomorphisms $s_{i,j}\:\hat g_{i}\sq\hat g_{j}\to
\hat g_{j}\sq\hat g_{i}$.  Use these in place of the commuting
isomorphisms in defining morphisms to products via the standard
algorithm.  The same proof shows the morphisms produced by the
algorithm are well-defined, since special properties of $\sigma$ were
not used.  The difference comes in associations.  As above, when
products are reduced in blocks specified by associations rather than
all at once, the ``commuting'' isomorphisms $s_{i,j}$ occur out of the
standard order.  Now, however, the crossing identity is no longer
valid so each of these out-of-order interchange contributes a
correction factor.  These are the new ingredients of the general case.

\subsubhead{2.5.7 Order conditions }\endsubsubhead The arguments of 2.5.6
give uniqueness, ie,  that there is a braided--monoidal equivalence
from a group--category $\Cal C$ to the standard one with the same
invariants.  However this implicitly uses the existence assertions,
that the invariants of $\Cal C$ satisfy the order conditions, and
conversely if a set of invariants satisfy the order conditions then
the twisted structure on $\Cal R[G]$ does in fact give a
braided--monoidal category.  We will discuss the cyclic case, ie,  the
$\sigma_{i}$ which commute a generator with itself, since this has the
extra factor of 2 and the connection to associativity.  The conditions
on $\sigma_{i,j}$ which commutes distinct generators are more routine
and are omitted.

Fix a generator of $G$, and drop the index $i$ from the notation
$g_{i}$.  Thus the generator is $g$, its order is $n$, $\hat g$ is the
chosen simple object in the equivalence class, $\Lambda\:1\to \hat
g^{n}$ is the chosen isomorphism implementing the order, and the
commutativity isomorphism $\hat g\sq \hat g\to \hat g\sq\hat g$ is
multiplication by $\sigma$.  Finally define $\alpha\in R$ so that the diagram
$$\CD \hat g @>\alpha >> \hat g\\
@VV {\Lambda\sq\id} V @VV{\id\sq\Lambda}V\\
(\hat g\sq\hat g^{n-1})\sq \hat g @>>>\hat g\sq(\hat g^{n-1}\sq \hat g)\endCD$$
commutes, where the top morphism is multiplication by $\alpha$ and the
bottom is the associativity isomorphism in the category.  The
conditions in 2.5.1 for a single generator are equivalent to:

\proclaim{Lemma} $\sigma^{n}=\alpha=\alpha^{-1}$\endproclaim
The hexagon axiom for commutativity isomorphisms asserts that the diagram
commutes (where unmarked arrows are associativities):
$$\CD
( g\sq g^{k-1})\sq g @>>>
 g\sq( g^{k-1}\sq g)\\
@VV{\sigma_{g,g^{k-1}}}V @VV{\sigma_{g,g^{k}}}V\\
( g^{k-1}\sq g)\sq  g @. ( g^{k-1}\sq g)\sq
 g\\
@VVV @VVV\\
 g^{k-1}\sq( g\sq  g)@>\sigma_{g,g}>> g^{k-1}\sq( g\sq  g)
\endCD$$
The reduction algorithm of 2.5.5 give canonical maps from a standard $
g^{k+1}$ into these objects, and we think of these as bases for $\hom(
g^{k+1},*)$.  If $k<n$ then the associativities are all ``identities''
(preserve these canonical bases).  The commuting maps multiply by
elements of $R$, so for these elements the diagram gives a relation
$\sigma_{g,g^{k}}=\sigma_{g,g^{k-1}}\sigma_{g,g}$.
$\sigma_{g,g}$ is multiplication by $\sigma$, so this subsitution and
 induction gives $\sigma_{g^{j},g^{k}}=\sigma^{jk}$, if $j, k <n$.

Now consider the diagram with $k=n$. The previous argument  still applies
to the left side and bottom, and shows the diagonal composition is
$\sigma^{n}$.  $\alpha$ is defined so the top associativity takes
$\Lambda\sq\id$ to $\alpha(\id\sq\Lambda)$.  We can evaluate the
upper right $\sigma$ term using the unit condition.  This condition
requires that the diagram commutes: $$\CD g\sq 1@>\sigma_{ g,1}>> 1\sq
g\\
@VV{{}\nu_{ g}}V @VV{{}_{ g}\nu}V\\
 g @>=>> g
\endCD$$
Composing the inverse of this with $ \Lambda\:1\to  g^{n}$ and using
naturality gives
$$\CD  g @>=>> g\\
@VV{\id\sq \Lambda}V @VV{ \Lambda\sq\id}V\\
   g\sq  g^{n} @>\sigma_{ g, g^{n}}>>  g^{n}\sq g\endCD$$
This shows the upper right side in the main diagram takes $\id\sq \Lambda$
to $ \Lambda \sq \id$. Therefore going across
the top and down the right side takes the standard generator to
$\alpha$ times the standard generator. Comparing with the other composition
gives $\sigma^n=\alpha$.

There is a second hexagon axiom in which $\sigma_{a,b}$ is replaced by
$\sigma_{b,a}^{-1}$.  The same argument applies to this diagram to give
$(\sigma^{-1})^n=\alpha$.  This completes the proof of the identity.
In fact this proof shows that the identities are exactly equivalent to
commutativity of the diagrams above, so the identity implies the
hexagon axioms.  To complete the argument it must be verified that the
formula for association in 2.5.2 satisfies the pentagon axiom if
$\alpha^{2}=1$.  This is straightforward so is omitted.

\subsubhead 2.5.8 Balance \endsubsubhead The final task is to show that
braided group--categories are  balanced, ie,
there is a functor $\tau(a)$ so that (writing the operations
multiplicatively)
$$\sigma(a,b)^{-1}=(\tau(a)\sq\tau(b))\sigma(b,a)\tau(ab)^{-1}.$$ In fact
$\tau(a)=\sigma(a,a)$ works.
\proclaim{ Lemma} The commuting
isomorphism $\sigma$ in a braided group--category satisfies
$$\sigma(ab,ab)=\sigma(a,a)\sigma(b,b)\sigma(a,b)\sigma(b,a)$$
\endproclaim
Note this relation would follow if $\sigma$ were bilinear, but this is
usually not the case.
\demo{Proof}In the following we use freely the fact that $R$ is a
commutative ring, so even though the identities are written
multiplicatively they can be reordered at will.  First, the hexagon
axiom for $(ab, a, b)$ gives
$$\sigma_{ab,ab} =
\sigma_{ab,a}\sigma_{ab,b}\alpha_{ab,a,b}^{-1}\alpha_{a,b,ab}^{-1}\tag{1}$$
Next the pentagon axiom for $(a,b,a,b)$ gives
$$\alpha_{a,ab,b}\alpha_{ab,a,b}^{-1}\alpha_{a,b,ab}^{-1}=\alpha_{a,b,a}^{-1}\alpha_{b,a,b}^{-1}.$$
Subsituting this into (1) gives
$$\sigma_{ab,ab}=(\sigma_{ab,a}\alpha_{a,b,a}^{-1})(\sigma_{ab,b}\alpha_{b,a,b}^
{-1})\tag{2}$$
In the inverse hexagon for $(a,b,a)$ two $\alpha$ terms cancel to give
$$\sigma_{a,ba}\alpha_{a,b,a}^{-1}=\sigma_{a,b}\sigma_{a,a}$$
Subsituting this, and the similar formula obtained by interchanging
$a$ and $b$, into (2) gives the identity of the lemma.\qed\enddemo

\head Homological field theories \endhead The ``theory'' based on $n^{th}$
homology is described in 3.1.  It is defined for general topological
spaces, but is not a field theory in this generality.  Criteria for
this are given in 3.1.3.  In particular the $H_{n}$ theory is modular on
$(n+1)$--complexes, but is a nonmodular field theory on $(n+2)$--manifolds.  In
3.2 the $H_{1}$  theory on 2--complexes is shown to agree with the categorical
construction using a group--category.  More general theories are
obtained in Section 4 by twisting the dual cohomology-based theories.

\subhead 3.1 The $H_{n}$ field theory \endsubhead The objective is to use
homology groups to define a topological field theory.  The definition
is given in 3.1.1, and hypotheses implying the field theory axioms are
given in 3.1.3.  Examples are given in 3.1.4, and in particular the
$H_{n}$ theory is a non-modular field theory on $M^{n+2}$ manifolds.
In 3.1.5 the $H_{1}$ theory on 2--complexes is shown to be the
category-based theory defined using the canonical group--category.  In
the following ``space'' will mean finite CW complex, ``subspace''
means subcomplex.  These assumptions imply that homology groups are
finitely generated, and pairs satisfy excision, long exact sequences,
etc.

\subsubhead 3.1.1 Definition \endsubsubhead
 Fix a commutative ring $R$,
a finite abelian group $G$ and a dimension $n$.  For a  pair
$(Y,W)$ define the ``state space'' by
$$Z(Y,W) = R[H_{n}(Y,W;G)].$$
Next suppose $X\supset Y_{0}\cup Y_{1}$ and $Y_{0}\cap Y_{1}=W$.
Then the induced homomorphism $Z_{X}\:Z(Y_{0},W)\to Z(Y_{1},W)$ is
defined by: for $y\in H_{n}(Y_{0},W;G)$,
$$Z_{X}(y)=\Sigma_{\{x\mid \partial_{0}x=-y\}}\partial_{1}x.$$
The $x$ in the sum are elements of $H_{n+1}(X,Y_{0}\cup Y_{1};G)$, and
the $\partial_{i}$ are boundary homomorphisms
$\partial_{i}\:H_{n+1}(X,Y_{0}\cup
Y_{n};G)\to H_{n}(Y_{i},W;G)$.

 $Z_{X}$ can be described a bit more explicitly using the exact sequence
 $$H_{n+1}(X)@>>>H_{n+1}(X,Y_{0}\cup Y_{1})@>\partial>>
 H_{n}(Y_{0})\oplus H_{n}(Y_{1})
  @>i>>H_{n}(X).$$
Let $k$ be the order of the image of $H_{n+1}(X)$ in $H_{n+1}(X,Y_{0}\cup
Y_{1})$.  Then
$$Z_{X}(y) = k\Sigma\{y_{1}\in H_{n}(Y_{1}) \mid i(y_{1}) = i(y)\}.$$

We want to find conditions under which this defines a topological
field theory, and when the theory is modular.

\subsubhead 3.1.2 Axioms \endsubsubhead
 Domain categories are defined in
\cite{Q} as the appropriate setting for topological field theories,
but full details are not needed here.  We take the objects (spacetimes) of the
category to be a subcategory $\Cal T$ of topological pairs $(X,Y)$.  The
boundary objects are the possible second elements $Y$.  The definition
above satisfies the tensor property (disjoint unions give tensor
products of state spaces, morphisms) on any $\Cal T$ because disjoint
unions give direct sums in homology.  The composition property
requires that if $X_{1}\:Y_{0}\to Y_{1}$ and $X_{2}\:Y_{1}\to Y_{2}$
are bordisms then $Z_{X_{2}}Z_{X_{1}} = Z_{X_{1}\cup X_{2}}$.  This is
not satisfied for completely general $\Cal T$.

 In a {\it modular\/} domain category three levels of objects are
specified. Boundary objects have
corner objects as their boundaries and certain identifications are
allowed.  A field theory on a modular domain category has relative
state spaces $Z(Y,W)$ defined for a (boundary, corner) pair, and
induced homomorphisms defined for boundaries with corners.  Here we
assume the extended boundary objects $(Y,W)$ are certain specified
topological pairs, glueing is the standard topological operation,
etc,  and then definition 3.1.1 is given in the modular formulation.
If $Z$ is a field theory on a modular domain category then for each
corner object $W$ the state space $Z(W\times I,W\times\partial I)$ has
a natural ring structure, and if $Y$ is a boundary object with
boundary $W_1\cup W_2$ then the state space $Z(Y,W_1\cup W_2)$ has
natural module structures over the corner algebras $Z(W_i\times
I,W_i\times\partial I)$.  A field theory is modular if the state space
of a glued object is obtained by ``algebraically'' glueing the state
space of the original object.  More specifically suppose $(Y,\partial
Y)$ is a boundary object with a decomposition of its boundary in the
corner category, $\partial Y = W_1\cup W_2\cup V$, and $W_1\simeq
\br{W_2}$.  Then there is a glueing in the category, $(\cup_WY,V)$,
and a natural homomorphism of state spaces
$$Z(Y,W_1\cup W_2\cup V)@>>> Z(\cup_WY,V).$$
The two copies of $W$ give two module structures on $Z(Y,W_1\cup W_2\cup
V)$  over the ring $Z(W\times
I,W\times \partial I)$, and the difference between the two vanishes in
$Z(\cup_WY,V)$.
This gives a factorization of the natural homomorphism through
$$\otimes_{Z(W\times
I,W\times \partial I)}Z(Y,W_1\cup W_2\cup V)@>>> Z(\cup_WY,V).\tag{$*$}$$
The field theory is said to be modular if this homomorphism is an isomorphism.

In the following $\Cal T$ is a domain category whose objects are
(certain specified) topological spaces.  Examples are given in 3.1.4.

\proclaim{3.1.3 Lemma} $Z$  satisfies the composition property (so defines
a field
theory) on
$\Cal T$ provided: if\/ $(X,Y)$ is a $\Cal T$  pair and\/
 $Y=Y_{1}\cup_{W}Y_{2}$ is a $\Cal T$ decomposition then
 $$H_{n+2}(X,Y_{1}\cup Y_{2};G)@>\partial>> H_{n+1}(Y_{1},W;G)$$
 is onto.  If $\Cal T$ is a modular topological domain category then $Z$ is
 modular provided in addition: if\/ $\cup_WY$ is a glueing in the
 boundary category, with boundary\/ $V$, then the homomorphism
$$H_{n+1}(\cup_WY,V;G)@>\partial>>H_n(W;G)$$
is onto.
 \endproclaim

\subsubhead 3.1.4 Examples\endsubsubhead
\roster \item $Z$ is a modular field theory on the modular domain category of
$(n+1)$--complexes, ie, with (objects, boundaries, corners) =
($(n+1)$--complexes,
$(n)$--complexes, $(n-1)$--complexes). Slightly more generally, it is
sufficient to have the homotopy type of complexes
of the indicated dimensions.  The composition and modularity conditions are
satisfied  because the groups involved are
all trivial.
\item $Z$ is a field theory on the domain category of oriented
$(n+2)$--manifolds, ie,  with (objects, boundaries) = ($(n+2)$--manifolds,
$(n+1)$--manifolds).  In this case $H_{n+2}(X,Y_{1}\cup Y_{2};G)$ and
$H_{n+1}(Y_{1},W;G)$ are both isomorphic to $G$ generated by the
respective fundamental classes, and the boundary homomorphism is an
isomorphism.  However the theory is not modular on the modular domain category
with corners $n$--manifolds.  The criterion given in the lemma fails because
$H_n(W;G)\simeq G$, and when $Y$ is obtained by identifying two copies
of $W$ the boundary homomorphism $\partial\:H_{n+1}(Y,\partial Y)\to
H_{n}(W)$ is trivial.  More directly, the theory is not modular
because the modularity construction does not account for the image of
the fundamental class of $W$ in $H_{n}(Y,\partial Y)$.  \endroster

 \demo{ Proof of 3.1.3} The composition property for
$(X_{1},Y_{0}\cup_{W}Y_{1})$
 and $(X_{2},Y_{1}\cup_{W}Y_{2})$ is that the functions
 $Z_{X_{1}\cup_{Y_{1}}X_{2}}$ and $Z_{X_{2}}Z_{X_{1}}$ agree. Both are
defined as sums of
$\partial_2$ of homology classes, so we need to show there is an
appropriate bijection between the
index sets.

 There is a commutative diagram with excision
 isomorphisms on the top and bottom, $$\CD
 H_{k}(X_{1}\cup_{Y_{1}}X_{2}, Y_{0}\cup Y_{1}\cup Y_{2})@>\simeq >>
 H_{k}(X_{1}, Y_{0}\cup Y_{1})\oplus H_{k}(X_{2},Y_{1}\cup Y_{2})\\
@VV{\partial}V @VV{\partial_{1}^{X_{1}}-\partial_{1}^{X_{2}}}V\\
H_{k-1}(Y_{0}\cup Y_{1}\cup Y_{2} , Y_{0}\cup Y_{2})@>\simeq >>
H_{k-1}(Y_{1},W)\endCD$$
Using this to replace terms in the long exact sequence of the triple
$X_{1}\cup_{Y_{1}}X_{2}\supset Y_{0}\cup Y_{1}\cup Y_{2} \supset
Y_{0}\cup Y_{2}$
gives $$\align H_{n+2}(X_{1}, Y_{0}\cup Y_{1})\oplus
H_{n+2}(X_{2},Y_{1}\cup Y_{2})& @>{\partial_1-\partial_1}>>H_{n+1}( Y_{1},
W) @>i>> \\H_{n+1}(X_{1}\cup_{Y_{1}}X_{2}, Y_{0}\cup Y_{2})
 @>j>> H_{n+1}(&X_{1}, Y_{0}\cup Y_{1})\oplus H_{n+1}(X_{2},Y_{1}\cup Y_{2})\\&
 @>{\partial_1-\partial_1}>>H_{n}( Y_{1}, W) \endalign$$
 The index set for the sum in $Z_{X_1\cup X_2}$ is the middle term, while
the index set for
the composition is the kernel of the lower boundary homomorphism. The
function $j$ between these is
onto by exactness.  For it to also
 be one-to-one we need $i=0$, or equivalently the upper boundary
 homomorphism is onto.  But this is the sum of two morphisms, both of
 which are onto by the hypothesis of the lemma, so it is onto.

Now consider the modular case. The ring structures are obtained by
applying $Z$ to $(W\times I)\times I$, regarded as a bordism rel ends
from $W\times I\cup W\times I$ to $W\times I$.  Similarly
$Y\cup_WW\times I\simeq Y$, so $Y\times I$ can be regarded as a
bordism $Y\cup W\times I\to Y$.  Applying $Z$ to this gives the module
structure.  In the case at hand $Z(W\times I,W\times \partial
I)=R[H_n(W\times I,W\times \partial I)]$, and the ring structure is
pointwise multiplication in the free module.  (This means if $v,w$ are
basis elements then $vw=0$ if $v\neq w$, and $vw=w$ if $v=w$.) There
are isomorphisms $H_n(W\times I,W\times \partial
I)@>\partial_i>>H_{n-1}(W)$ for $i=0,1$, and $\partial_0=-\partial_1$.
The module structure on $R[H_n(Y,V\cup W)]$ using the 1 end of
$W\times I$ is: if $y\in H_n(Y,V\cup W)$, $v\in H_{n-1}(W)$ then $vy=0$ if
$\partial_Wy\neq v$, and $vy=y$ if $\partial_Wy = v$.  Using the other
end of $W\times I$ gives 0 or $y$ depending on whether or not
$\partial_Wy = -v$.

This description of the ring and module structures identifies the algebraic
glueing on the left in
the modularity criterion ($*$) as the free module generated by $y\in
H_n(Y,W_1\cup
W_2\cup V;G)$ satisfying $\partial_{W_1}y=-\partial_{W_2}y$.

Now consider the long exact sequence of the triple $\cup_WY\supset V\cup
W\supset V$:
$$\align H_{n+1}(\cup_WY,V\cup
W)@>\partial>>H_n(W)@>>>H_n(\cup_WY,V)@>>>&H_n(\cup_WY,V\cup
W)\\&@>\partial>>H_{n-1}(W)\endalign$$
The state space of the geometric glueing is generated by the third term,
while we have identified
the algebraic glueing as generated by the kernel of $\partial$ in the
fourth term. The homomorphism
of ($*$) is induced by the set-level inverse of the third homomorphism, so
we need to show the third
homomorphism is an isomorphism onto the kernel of $\partial$. Exactness
implies it is onto. For
injectivity we need the second homomorphism to be 0, or equivalently the
first $\partial$ to be
onto. But this is exactly the hypothesis of the lemma.
\qed\enddemo

\subhead 3.2 Connections to categories\endsubhead This gives the first
direct connection between the homological theories and categorical
constructions.  The general case is in Section~4.

 \proclaim{Proposition} The
canonical untwisted group--categories are the only ones that define
modular field theories on 2--complexes, and the corresponding field
theories are the $H_{1}$ theories of 3.1.1.  \endproclaim

\demo{Proof} The categorical input for fields on 2--complexes is a
symmetric monoidal category satisfying a symmetry condition.
Symmetric monoidal group--categ\-or\-ies are classified in 2.3(3), or
2.5.1(3).  The first part of 3.1.4 corresponds to the fact that of
these only the canonical examples satisfy the symmetry conditions.

The symmetry condition concerns nondegenerate pairings.  A
nondegenerate pairing on $a$ is another object $\br a$ and morphisms
$$\align \Lambda_{a}\:&1\to \br{a}\sq a\\
\lambda_{a} \:&a\sq \br{a}\to 1\endalign$$
satisfying
$$ a\simeq a\sq 1@>\id\sq\Lambda_{a}>>a\sq(\br a\sq
a)@>\text{associate}>>(a\sq\br a)\sq a@>\lambda_{a}\sq \id>>1\sq a\simeq a $$
$$\br a\simeq 1\sq a@>\Lambda_{a}\sq\id>>(\br a\sq a)\sq \br
a@>\text{associate}>>\br a\sq(a\sq \br a)@>\id\sq\lambda_{a}>>\br a\sq
1\simeq \br a$$
are both identity maps.  The construction requires a fixed choice of
pairings on the simple objects.  This is equivalent to an additive
assignment of pairings to all objects, and this in turn is equivalent to
a ``duality'' functor making the category ``autonomous'', \cite{S} or
a nondegenerate trace function.

The construction requires the symmetry condition  $\lambda_{\br
a}=\lambda_{ a}\sigma_{ \br a,a}$.
If $a\neq \br a$ then we can arrange this to hold by taking it as the
definition of $\lambda_{\br a}$.  If $a=\br a$ the condition is
equivalent to $\sigma_{a,a}$ being the identity.  But this is the only
possibly nontrivial invariant in 2.5.1(3), so the category is standard.

Now we show that the $H_{1}$ theory corresponds to the standard
group--category.  One way to do this is to go through the construction
[\Cite{Q}, \Cite{B}] and see homology emerge.
 This is illuminating but too
long to reproduce here.  Instead we use the reverse
construction, extracting a category from a field theory.
This goes as follows: let ``$\pt$'' denote the connected corner object
in the domain category.  The state space of $Z(\pt\times I)$ has a
natural ring structure, and additively the category is the category of
modules over this ring.  The state space of the cone on three points
$Z(c(3))$ has three module structures over the ring.  The product on
the category is defined by tensoring with this trimodule.

$Z(\pt\times I) = R[H_{1}(I,\partial I; G)] = R[G]$.  The ring structure is
obtained by considering the boundary of $I^{2}$ as the union of three
intervals, with two incoming and one outgoing.  If $(g,h)\in
H_{1}(I,\partial I)\oplus H_{1}(I,\partial I)$ then the image in
$R[H_{1}(I,\partial I)]$ is obtained by summing over elements of
$H_{2}(I\times I,\partial I\times I;G)$ whose restrictions to the
incoming boundary intervals is $g$ and $h$.  $H_{2}(I\times I,\partial
I\times I;G) = G$, and the restrictions are identities.  Thus $(g,h)$
goes to 0 if $g\neq h$, and to $g$ if they are the same.  The ring is
therefore $R[G]$ with componentwise multiplication.  There is an
anti{}involution on this
ring induced by interchanging ends of the interval.  This is the
involution on $R[G]$ induced by inverse in $G$.  The category of
modules over this ring is exactly the $G$--graded (left) $R$--modules.
Denote the category by $\Cal C$.  Simple objects are $R[g]$ as in 2.2.1:
a copy of $R$ on which multiplication by $h\in G$ is zero if $h\neq g$
and is the identity if $h=g$.

Now let $c(3)$ denote the cone on three points.  The standard cell
structure is three invervals joined at a point.  Using cellular chains
gives an explicit description of $H_{1}(c(3),3;G)$ as $\{(a,b,c)\in
G^{3}\mid abc=1\}$ (the three generators correspond to the three
1--cells, the relation comes from the boundary homomorphism to the
chains on the vertex).  The three (left) module structures over
$R[H_{1}(I,\partial I; G)]$ are defined by glueing intervals on the
three endpoints.  Thus in the first structure $g$ in the ring takes
$(a,b,c)$ to 0 if $g\neq a$, and $(a,b,c)$ if $g=a$.  The product on
the category
$$\Cal C\times \Cal C @>>>\Cal C$$
is defined by: begin with $M$ and $N$ left modules over the ring.
Convert these to right modules using the antiinvolution in the ring,
and tensor with the first two module structures on $Z(c(3))$.  Then
$M\sq N$ is the result, with respect to the third module structure.
Now we can work out the product of two simple objects $R[g]\sq R[h]$.
The involution converts these to right modules on which $g^{-1}$,
$h^{-1}$ respectively act nontrivially.  Tensoring with the first two
coordinates in $R[\{(a,b,c)\in G^{3}\mid abc=1\}]$ kills everything
with $a\neq g^{-1}$, $b\neq h^{-1}$, so leaves exactly $R[gh]$.
Therefore the product is the standard product in $\Cal R[G]$.

This does not yet identify the category as standard: according to
2.3.2 any group--category is equivalent to the standard one with the
standard product.  The differences are in the associativity and
commutativity structures.  Here commutativity comes from the
involution on the cone on three points that interchanges the two
``incoming'' ends.  This interchanges two of the 1--cells in the cell
structure, so interchanges the corresponding generators in the
cellular 1--chains.  Thus in homology it interchanges the first two
coordinates in $\{(a,b,c)\in G^{3}\mid abc=1\}$.  Following through
the tensor product  gives the standard ``trivial'' commuting
isomorphism for $R[g]\sq R[h] = R[g]\otimes R[h]$.  This finishes the
argument because the commutativity
determines the associativity.  Standardness of associativity is also
easy to see directly: associating isomorphisms come from two ways to glue
together two cones on three points to get (up to homotopy) the cone on
four points.  Following this through gives the standard trivial
associations.\qed\enddemo

\head Cohomological field theories \endhead
Homology will probably be the most natural setting for field theories, but
so far only the fields for  standard
group--categories can be described this way. In this section we restrict to
manifolds and show how to twist the dual
theory theory based on cohomology. More specifically  fix a space  $E$ with
two nonvanishing homotopy
groups
$\pi_d E=G$ and
$\pi_{n+d}E=\units(R)$, and suppose $E$ is simple if $d=1$. We construct
state spaces and induced
homomorphisms from homotopy classes of maps to this space. A simple case is
described in 4.1.1 to show this gives a
twisted version of the $H_n(*;G)$ theory. The full definition occupies the
rest of 4.1. The field axioms and
modularity are verified in 4.2.   The $n=1$ cases are shown to be
Reshetikhin--Turaev constructions from
group--categories in 4.3.

\subhead 4.1 The definition \endsubhead
The general construction is a bit complicated so we begin with a special
case in 4.1.1. The domain category for
the theory is defined in 4.1.2; the special case of 4.1.1 is supposed to
explain why this is the right choice. Once
the objects are known the full definition can be presented.

\subsubhead 4.1.1 A special case\endsubsubhead
The Postnikov
decomposition for the fixed space $E$ is
$$B^{n+d}_{\units(R)}@>>> E@>>> B^d_G @>k>>B^{n+d+1}_{\units(R)}.$$
The first space $B^{n+d}_{\units(R)}$ has the structure of a topological
abelian group, and the last space is the
classifying space for principal bundles with this group. In particular $E$
is a principal bundle with an action of
$B^{n+d}_{\units(R)}$.

Now suppose
$Y$ is a connected oriented manifold of dimension $n+d$. The
 group $[Y/\partial Y, B^{n+d}_{\units(R)}] = H^{n+d}(Y,\partial
Y;\units(R))$  is dual to
$H_0(Y;\units(R)) = \units(R)$. This acts on the set of homotopy classes
$[Y/\partial Y,E]$ and
 the quotient of this action is (when $G$ is abelian)
$$[Y/\partial Y,B^d_G]=H^d(Y,\partial y;G) \simeq H_n(Y;G).$$ Define the
state space
$Z(Y)$ to be the set of functions
$[Y/\partial Y,E]\to R$ that commute with the action of $\units(R)$.

If $E$ is the product $B^{n+d}_{\units(R)}\times B^d_G$ then the homotopy
classes are also a product $[Y/\partial
Y,E] = \units(R)\times  H_n(Y;G)$ and the set of $\units(R)$--maps is
$R[H_n(Y;G)]$, exactly the definition
of Section 3. Thus the
$k$--invariant of $E$ gives a way to twist the $R$--module generated by
$H_n(Y;G)$. In the present case ($Y$ connected)
this can also be described as: $[Y/\partial Y,E]$ is a principal
$\units(R)$ bundle over
$H_n(Y;G)$. The state space is the space of sections of the associated
$R$--bundle.

 Note to get the key canonical
identification of $[Y/\partial Y, B^{n+d}_{\units(R)}]$ with $\units(R)$ we
needed the boundary
objects to be oriented manifolds of dimension $n+d$.

\subsubhead 4.1.2 The domain category\endsubsubhead
The  field theory will be defined on $(n+1+d)$--dimensional thickenings of
$(n+1)$--complexes. The
definitions of state spaces and induced homomorphisms use only the manifold
structure. Restrictions on the homotopy
dimension are needed for the field axioms to be satisfied.
\roster
\item Corner objects are compact oriented $(n+d-1)$--manifolds with the
homotopy type of an
$(n-1)$--complex, together with a set of maps
$w_i\:W/\partial W\to E$, one in each homotopy class;
\item relative boundary objects are compact oriented $(n+d)$--manifolds with
the homotopy type of an $n$--complex, with
boundary given as a union
$\partial Y=\hat\partial Y\cup W$ of submanifolds, and $W$ has the
structure of a corner object (ie, homotopy
dimension $n-1$ and a choice of maps $w_i$); and
\item ``spacetime'' objects are compact oriented $(n+d+1)$--manifolds with
homotopy type of $(n+1)$--complexes,
boundary given as a union
$\partial X=\hat\partial X\cup Y$ of submanifolds, and $Y$ having the
homotopy type of an $n$--complex.
\endroster

The ``internal'' boundary (in the domain category) of an object
$(X,\hat\partial X\cup Y)$ is $Y$ with $\hat\partial
Y=\partial Y$ and $W=\emptyset$. The internal boundary of $Y$ with
$\partial Y=\hat\partial Y\cup W$ is $W$.
Morphisms are orientation-preserving homeomorphisms, required to commute
with the fixed reference maps on corners.
The choices of maps in (1) are typical of the rigidity seen in corner
objects, see
\cite{Q}. The involution $X\mapsto \overline X$ is defined by reversing the
orientation.

\subsubhead 4.1.3 The definition\endsubsubhead
Suppose $Y$ with $\partial Y= \hat\partial Y\cup W$ and $w_i\:W/\partial
W\to E$ is a relative boundary object.
Define $[Y/\hat\partial Y,E]_0$ to be maps that agree with one of the
standard choices on $W$, modulo homotopy
rel $\partial Y$. The group $[Y/\partial
Y,B^{n+d}_{\units(R)}]=H^{n+d}(Y,\partial Y;\units(R))$ acts on this set,
as in 4.1.1. {\it Caution\/}: the group operation in $\units(R)$ is written
{\it multiplicatively\/}. The operations in cohomology groups and their
action on homotopy classes into $E$ are
therefore also written multiplicatively. Define
$$\epsilon\: H^{n+d}(Y,\partial Y;\units(R))\to \units(R)$$
 by evaluation on the fundamental class
of $Y$. When $Y$ is connected (as in 4.1.1) this is an isomorphism, but we
do not assume that here. Define the state
space for the theory by
$$Z(Y,W)=\hom_{\epsilon}([Y/\hat\partial Y,E]_0,R)$$
where $\hom_{\epsilon}$ indicates functions $\alpha\:[Y/\hat\partial
Y,E]_0\to R$ so that if
 $f\in [Y/\hat\partial Y,E]_0$ and $a\in H^{n+d}(Y,\partial
Y;\units(R))$ then $\alpha(a f) = \epsilon(a) \alpha(f)$.

Now we define induced homomorphisms. The general modular setting is an
object with boundary divided into ``incoming''
and ``outgoing'' pieces, and the incoming boundary further subdivided.
Specifically suppose $Y_1$ is a relative
boundary object with corner a disjoint union $W_1\sqcup W'_1\sqcup W_2$, an
isomorphism $W'_1\simeq \overline W_1$ is
given, and $\cup_{W_1}Y_1$ is the object obtained by identifying $W'_1$ and
$W_1$. Suppose $Y_2$ is a boundary object
with corner $W_2$, and finally $X$ is an object with internal boundary
$(\cup_{W_1}\overline Y_1)\cup_{W_2}Y_2$. Then
we define
$$Z_X\:Z(Y_1,W_1\sqcup W'_1\sqcup W_2)@>>> Z(Y_2, W_2)$$
as follows. An element in the domain is a function
$\alpha\:[Y_1/\hat\partial Y_1, E]_0\to R$. The output is
a function $[Y_2/\hat\partial Y_2, E]_0\to R$, so we can define it by
specifying its value on a map
$f\:Y_2/\hat\partial Y_2\to E$. We first suppose each component of $X$
intersects either $Y_1$ or $Y_2$. Then
$$Z_X(\alpha)(f)=\Sigma_{[g]}\epsilon_2(a)\alpha(\hat g|Y_1).$$
The sum is over homotopy classes of $g\:X/\hat\partial X\to B^d_G$ whose
restriction to $Y_2$ is homotopic to the
projection of $f$. When
$G$ is abelian (eg if
$d>1$) this is
dual to index set used in the homological version.
$\hat g\:X/\hat\partial X\to E$ is a lift of $g$ which is standard on $W_1$
and $W_2$, and $a\in H^{n+d}(Y_2,\partial
Y_2;\units(R))$ so that
$a\cdot
\hat g|Y_2\sim f$.  When each component of $X$ intersects either $Y_1$ or
$Y_2$ such a lift
exists, and since $\hat g|Y_2$ and $f$ project to homotopic maps in $B^d_G$ they differ by the action of some such
element $a$.   

If $Y_1$ and $Y_2$ are empty then  we define an element of $R$ by
$$\hat Z_X=\Sigma_{[g]}k(g)([X]).$$
Here the sum is again over $X/\partial X\to B^d_G$, $k\:B^d_G\to
B^{n+d+1}_{\units(R)}$ is the $k$--invariant of $E$,
 (see 2.3.2 and 4.1.1) and $k(g)([X])$ is the evaluation of the resulting
cohomology class on the fundamental class
of $X$.

Now define $Z_X$ for general $X$. Write $X$ as $X_1\sqcup X_2$, where $X_1$
are the components intersecting $Y_1\cup
Y_2$ and $X_2$ are the others. If $X_1$ is nonempty define $Z_X$ as
$Z_{X_1}$ multiplied by $\hat Z_{X_2}$. If
$X_1$ is empty then $Z(Y_i)$ are canonically identified with $R$ and $Z_X$
is multiplication by $\hat Z_X$.

\proclaim{4.1.4 Lemma}  $Z_X$ is well-defined, and takes values in
$Z(Y_2,W_2)$. \endproclaim

\demo{Proof} The things to be checked are that $\epsilon_2(a)\alpha(\hat
g|Y_1)$ does not depend on the choice of
lift $\hat g$ and $a$, and that the resulting function $[Y_2/\hat\partial
Y_2,E]_0\to R$ commutes
appropriately with the action of $\units(R)$.

Suppose $g'$ is another lift of a map $g$. There is $b\in
H^{n+d}(X,\hat\partial X;\units(R))$ with $g'=b\cdot \hat
g$.
Denote the restrictions of $b$ to $Y_1$ and $Y_2$ by $b_1$ and $b_2$
respectively, then we have $f\sim a\cdot \hat
g|Y_2\sim a(b_2)^{-1}(b_2)\cdot(\hat g|Y_2)\sim a(b_2)^{-1}\cdot(b\cdot
\hat g)|Y_2\sim a(b_2)^{-1}\cdot( g')|Y_2$.
Therefore the element of $H^{n+d}(Y_2,\partial Y_2;\units(R))$ associated
to $g'$ is $a b_2^{-1}$, and the
corresponding contribution to $Z_X$ is $\epsilon_2(a
b_2^{-1})\alpha(g'|Y_1)$.  Since $\epsilon_2$ is a homomorphism
and $\alpha$ is an $\epsilon_1$--homomorphism,
$$\epsilon_2(ab_2^{-1})\alpha((b\cdot \hat
g)|Y_1)=\epsilon_2(a)\epsilon_2(b_2^{-1})\epsilon_1(b_1)\alpha(\hat
g|Y_1).$$ Thus we have to show $\epsilon_2(b_2)^{-1}\epsilon_1(b_1)=1$.
$\epsilon$ is defined by evaluation on
fundamental classes. The orientation of $Y_2$ is the opposite of the
induced orientation of $\partial X$, and the
complement of $Y_1\cup Y_2$ in $\partial X$ is taken to  the basepoint. Thus
$\epsilon_2(b|Y_2)^{-1}\epsilon_1(b|Y_1)$ is obtained by evaluating
$b|\partial X$ on the fundamental class of
$\partial X$. But $b|\partial X$ extends to a cohomology class ($b$) on
$X$, and the image of $[\partial X]$ in the
homology of $X$ is trivial (it is the boundary of the fundamental class of
$X$). Thus the evaluation is trivial; 1
since we are writing the structure multiplicatively.

To complete the lemma we show $Z_X(\alpha)$ is an $\epsilon_2$--morphism.
Suppose $f$, $\alpha$ as above, and $c\in
H^{n+d}(Y_2,\partial Y_2;\units(R))$. Then
$$\align
\epsilon_2(c)Z_X(\alpha)(f)&=\epsilon_2(c)\Sigma_{[g]}\epsilon_2(a)\alpha(\hat g
|Y_1) \\
&=\Sigma_{[g]}\epsilon_2(ac)\alpha(\hat g|Y_1)\\
&= Z_X(\alpha)(c\cdot f)
\endalign$$
The third line is justified by the fact that $a\cdot \hat g|y \simeq f$ if
and only if
$ac\cdot
\hat g|y \simeq c\cdot f$
\qed\enddemo

\subhead 4.2 The field axioms\endsubhead
We will not be as precise as in 3.1.3 about the exact conditions for field
axioms, but concentrate on the case of
interest. We continue the standard assumption that $E$ has two nontrivial
homotopy groups, $G$ in dimension $d$ and
$\units(R)$ in dimension $n+d$.
\proclaim{4.2.1 Proposition}  $Z$ defined in 4.2 is a modular field theory
on $(n+d+1)$--dimensional thickenings
of $(n+1)$--complexes. If $E$ is a product (and $G$ abelian if $d=1$) then
$Z$ is equal to the homological theory
of~3.1.
\endproclaim
\demo{Proof}
Consider the composition of $X_1\:Y_1\to Y_2$ and $X_2\:Y_2\to Y_3$.
Suppose first that each component of
$X_1\cup_{Y_2}X_2$ intersects either $Y_1$ or $Y_3$. In this case
$Z_{X_1\cup X_2}$ and $Z_{X_2}Z_{X_1}$ are given
by sums over  $[(X_1\cup X_2)/\hat\partial, B^d_G]$ and $[X_1/\hat\partial,
B^d_G]\times [X_2/\hat\partial, B^d_G]$
respectively. Since these are dual to the index sets used in the
homological theory, that proof shows that under the
given dimension restrictions the natural function between the two is a
bijection.Thus we need only show that the
corresponding terms in the sum are equal. Suppose $\alpha\in Z(Y_1)$ and
$f\in [Y_3/\hat\partial Y_3,E]_0$, and
consider the image of $\alpha$ evaluated on $f$.  Choose an element
$g\:(X_1\cup X_2)/\hat\partial\to B^d_G$ in the index set, and let
$\hat g$ be a lift, with $a\in H^{n+d}(Y_3,\partial Y_3;\units(R))$ so that
$a\cdot (\hat g|Y_3)\sim f$. The term in
$Z_{X_1\cup X_2}$ is $\epsilon_3(a)\alpha(\hat g|Y_1)$. Use restrictions of
$\hat g$ as lifts of the restrictions
of $g$ to $X_1$ and $X_2$. Since these agree on $Y_2$ there is no
$\epsilon_2$ correction factor, and the
corresponding term in
$Z_{X_2}Z_{X_1}$ is exactly the same.

Now consider a component of $X_1\cup X_2$ disjoint from $Y_1$ and $Y_3$, so
we want to show that $Z_{X_2}Z_{X_1}$ is
multiplication by the ring element $\hat Z_{X_1\cup X_2}$. If the union is
disjoint from $Y_2$ as well then it lies
entirely in one piece and this is trivially true. Thus suppose
$X_1\:\emptyset\to Y_2$, $X_2\:Y_2\to \emptyset$, and
$Y_2$ intersects each component of the union. Again the index sets match up
so we show corresponding terms are
equal. Choose a map $g\:(X_1\cup X_2)/\partial(X_1\cup X_2)\to B^d_G$, and
choose lifts $\hat g_1$ and $\hat g_2$ of
the restrictions to the two pieces. Note $g$ itself may not lift, so the
two lifts may not agree on $Y_2$. Let $a$
be a class with $a\cdot (\hat g_1|Y_2) = \hat g_2|Y_2$. Then we want to
show $a$ evaluated on the fundamental class
of $Y_2$ is the same as $kg$ evaluated on the fundamental class of $X_1\cup
X_2$.

 For convenience insert a collar on $Y_2$, so the union is $X_1\cup
Y_2\times I\cup X_2$. Now consider the lifts on the pieces as a lift of $g$
on the disjoint union. This lift gives a
factorization of $kg$ through $Y_2\times I/\partial$:
$$\CD (X_1\sqcup X_2)/\hat\partial @>{\hat g_1\sqcup \hat g_2}>>E\\
@VVV @VVV\\
(X_1\cup Y_2\times I\cup X_2)/\partial @>g>>B^d_G\\
@VVV@VV{k}V\\
Y_2\times I/\partial@>>>B^{n+d+1}_{\units(R)}\endCD$$
The lower map gives $kg$ as the image of an element in $H^{n+d+1}(Y_2\times
I,\partial(Y_2\times I;\units(R))$. The
suspension isomorphism
$$H^{n+d}(Y_2,\partial Y_2;\units(R)) \to H^{n+d+1}(Y_2\times
I,\partial(Y_2\times
I;\units(R))$$
 takes $a$ to this element. To see this, interpret the first class as the
classifying map for a
principal bundle over
$Y_2\times I$ filling in between the restrictions of $\hat g_i$ to $Y_2$.
The homotopy extension property for
principal bundles shows this is the mapping cylinder of a bundle
isomorphism, which must be the one classified by
$a$.  Since evaluation of $a$ on $[Y_2]$ is equal to the evaluation of the
suspension of $a$ on $[Y_2\times I]$, it
follows that $\epsilon_2(a)=kg([X])$.

This completes the proof of the composition property for induced
homomorphisms.
The proof of modularity is similar to the homology case, and in fact the
algebra associated to a corner object is
exactly the same.

Suppose $W$ is a corner object, so an oriented $(n+d-1)$--manifold with the
homotopy type of a $(d-1)$--complex
 and chosen representatives $w_i$ for homotopy classes $[W/\partial W, E]$.
The first claim is that there is a
canonical isomorphism
$$Z(W\times I)@>{\sim}>>\text{functions}([W/\partial W,B_G^d],R),$$   and
this takes the corner algebra
structure to the product induced by multiplication in $R$.  The definition
of $Z(W\times I)$ is
$\hom_{\epsilon}([(W/\partial W)\times I,B^d_G]_0$, where the subscript $0$
indicates that the restrictions to
$W\times \{0,1\}$ are images of  standard representatives $w_i$. The first
point is that the homotopy
$(d-1)$--dimensionality of $W$ implies that
$[W/\partial W,E]\to [W/\partial W,B^d_G]$ is a bijection. Since the
restrictions to the ends of a map $(W/\partial
W)\times I\to B^d_G$ are homotopic, this means the maps on the ends are
actually equal. Next, again using
dimensionality, a map
$(W/\partial W)\times I\to B^d_G$ which is equal to $pw_i$ on each end is
itself homotopic rel ends to the map which
is constant in the $I$ coordinate. This map has a canonical lift to $(W/\partial W)\times I\to E$ which is standard
on the ends, namely $w_i$ applied to
projection to the first coordinate. Applying $H^{n+d}(W\times
I,\partial(W\times
I),\units(R)$ to this gives a surjection  $[W/\partial W,B^d_G]\times
H^{n+d}(W\times I,\partial(W\times
I),\units(R)\to [(W/\partial w)\times I,E]_0$. Applying
$\hom_{\epsilon}(*,R)$ to this gives the required bijection.

The algebra structure in the algebra, or more generally the action on a
state space, is described as follows: Suppose
$Y$ is a relative boundary object with $\partial Y=\hat\partial Y\cup
W_1\cup W_2$. A function $\tau\:[W_1/\partial
W_1,B^d_G]\to R$ acts on  $\alpha\:[Y/\hat\partial Y,E]_0$\break $\to R$ to give
another function like $\alpha$. The
new function can be specified by its action on $f\in [Y/\hat\partial
Y,E]_0$, by
$$(\tau\cdot \alpha)(f)=\tau(f|W_1)\alpha(f).$$

Finally we prove modularity. Suppose $Y$ has corners $W\sqcup  \overline
W\sqcup W_2$, and let $\cup_WY$ be the
boundary object obtained by identifying the copies of $W$. The homomorphism
of state spaces induced by this glueing is
$$\hom_\epsilon([Y/\hat\partial Y,E]_0,R)@>>>
\hom_\epsilon([\cup_WY/\hat\partial Y,E]_0,R).$$
This is induced by a ``splitting'' function
$$[\cup_WY/\hat\partial Y,E]_0@>>> [Y/\hat\partial Y,E]_0$$
defined as follows. Suppose $f\:\cup_WY/\hat\partial Y\to E$ is standard on
$W_2$. The restriction to $W$ is
homotopic to a standard map. Use this to make $f$ standard on $W$, then
split along $W$ to obtain
$f'\:Y/\hat\partial Y\to E$ standard on $W\sqcup \overline W\sqcup W_2$.
The dimensionality hypotheses can be used as
above to show $f'$ is well-defined up to homotopy rel boundary, and the
splitting function is a bijection onto the
subset of
$g\:Y/\hat\partial Y\to E$ satisfying $g|W=g|(\overline W)$. Therefore to
show the algebraic glueing map
$$\otimes_{Z(W\times I)}Z(Y)@>>>Z(\cup_WY)$$
is an isomorphism we need to show that dividing by the difference between
the two $Z(W\times I)$--module structures
divides out exactly the functions supported on the complement of the image
of the splitting function. These
functions are sums of ``delta'' functions: suppose $g$ has $g|W\neq
g|(\overline W)$. Define $\alpha_g$ to take $g$ to
1, extend to an $\epsilon$--morphism on $H^(Y/\partial Y;\units(R))\cdot g$,
and define it to be 0 elsewhere. It is
sufficient to show these functions get divided out. Dividing by the
difference between the module structures divides
all elements of the form $f\mapsto (\tau(f|W) -\tau(f|\overline
W))\alpha(f)$. For the particular $g$ under
consideration there is a function $\tau$ with $\tau(g|W)=1$ and
$\tau(g|\overline W)=0$. Using this $\tau$ and the
delta function $\alpha_g$ gives
$$f\mapsto \cases
\alpha_g(f)&\text{ if $f$ is a multiple of $g$, since }\tau(f|W)
-\tau(f|\overline W)=1\\
 0 &\text{ otherwise , since }\alpha_g(f)=0
\endcases$$
But this is exactly $\alpha_g$, so $\alpha_g$ is divided out.
\qed\enddemo

\subhead 4.3 Relations to group--categories\endsubhead
Suppose $E$ is a space as above with $n=1$, so $\pi_d(E)=G$ and
$\pi_{d+1}(E)=\units(R)$.
In Section 2.3 these spaces are shown to correspond to group--categories
with various degrees of commutativity. The cohomological construction of
Proposition 4.2.1 gives a modular field theory on the domain category whose
objects, boundaries, corners are manifolds of dimension $(d+2,d+1,d)$ and
homotopy type of complexes of dimension
$(2,1,0)$ respectively.
Specifically we have:
\medskip
\hfil\vbox{\offinterlineskip
\halign{\vrule#&\strut\quad\hfil#\quad&\vrule#&\quad\hfil#\hfil
\quad&\vrule#&\quad\hfil#\hfil\quad&\vrule#\cr &$d$&&category structure
&&fields on & \cr
\noalign{\hrule}
&1&&associative &&$(3,2,1)$--thickenings &\cr
&2&&braided--commutative &&$(4,3,2)$--thickenings &\cr
&$\geq 3$&&symmetric &&$(d+2,d+1,d)$--thickenings &\cr
\noalign{\hrule}\cr
}}\hfil
\medskip
\noindent On the category side the independence of $d$ when $d\geq 3$ comes
from stability of group cohomology under
suspension. For fields,  cartesian product with $I$ gives a ``suspension''
functor of domain categories, from
$d$--thickenings to $(d+1)$--thickenings. Composition with this gives a
suspension function on field theories, from
ones on
$(d+1)$--thickenings to ones on
$d$--thickenings.  When
$d\geq3$, suspension is an equivalence of domain categories (ie, all
thickenings are isomorphic to products in an
appropriately canonical way), so  it induces a bijection of field theories.

There is also a Reshetikhin--Turaev type construction that uses a category
to define a field theory
on the same objects. Here we show that the two field theories agree.

\proclaim{Proposition} Suppose $\Cal G$ is a group--category corresponding
to a space $E$. The cohomological
field theory defined in 4.2 using $E$ is the same as the Reshetikhin--Turaev
theory defined using $\Cal G$.
\endproclaim

\demo{Proof} We will not prove this directly, but rather use the fact (as
in 3.2) that the
category can be recovered (up to equivalence) from the field theory.
Specifically the category $\Cal G$ is additively
equivalent to the category of representations of the corner algebra of a
thickening of a point, with product
structure induced by the state space of thickenings of the cone on three
points.

Fix a connected corner object: a copy of $D^d$ with specific choices of
representatives
$\hat g\:D^d/\partial D^d\to E$ for each homotopy class $g\in
[D^d/\partial D^d, E]=\pi_d(E)=G$.
The  algebra structure on
$Z(D^d\times I,D^d\times\{0,1\})$ is identified in the proof of 4.2.1 as
the set of functions $G\to R$, with product
given by product in $R$, or alternatively $R[G]$ with componentwise
multiplication.
Representations of this are exactly $\Cal R[G]$, so this gives an additive
equivalence $\Cal G\to \Cal R[G]$.

To determine the product structure we choose data as in 2.3.3--4: for each
pair $g,h$ choose a homotopy $m_{g,h}\:\hat
g\hat h\simeq
\widehat{gh}$. The left side of this expression is the product of homotopy
classes in $\pi_d(E)$, while the
right side is the given representative of the product in $G$. Let $Y$
denote the thickening of the cone on three
points, so $Y\simeq D^{d+1}$ with internal boundary $3D^d\subset \partial
D^{d+1}$ and $\hat \partial Y$ the
complement. The next object is to describe $Z(Y)$ with its three module
structures over the corner algebra. Reverse
the orientation on two of the boundary components (to switch the module
structure from left to right). A map
$Y/\hat\partial Y\to E$ that restricts to $\hat g\sqcup \hat h$ on the
incoming boundaries of $Y$ gives a homotopy to
the restriction to the third component. This identifies the third
restriction as $\widehat{(gh)^{-1}}$. The inverse
comes from the fact that all the components of $\partial Y$ have the
induced orientation, while in $D^d\times I$ one
of the ends has the reverse orientation. We have specified one such map, namely
$m_{g,h}$, and all others with this restriction are obtained (up to
homotopy) by the action of
$H^{d+1}(Y,\partial Y;\units(R))=\units(R)$. Thus the choices $m_{gh}$ give
a bijection
$$[Y/\hat\partial Y,E]_0\sim \cup_{g,h}\units(R).$$
The state space $Z(Y)$ is the set  $\hom_{\epsilon}([Y/\hat\partial
Y,E]_0,R)$, so the bijection gives an
identification $Z(Y)=R[G\times G]$. The three (left) module structures are:
on a summand $R[(g,h)]$, $f\in G$ acts by
the delta function $\delta_{f,g}$, $\delta_{f,h}$, and
$\delta_{f,(gh)^{-1}}$. Switch the first two to right
structures by reversing the orientation, and replace $g, h$ by $g^{-1},
h^{-1}$. This gives an identification in which
the right structures on $R[(g,h)]$ are $\delta_{f,g}$ and $\delta_{f,h}$
respectively, and the left structure is
$\delta_{f,gh}$. Now suppose $a_g$ and $a_h$ are simple modules in $\Cal
R[G]$. Their $Z$--product is
$Z(Y)\otimes_{Z(W\times I)^2}(a_g\otimes a_h)$. The description of $Z(Y)$
shows this is a free based module of rank 1,
canonically isomorphic to $a_{gh}$. This gives a natural isomorphism
between the $Z$ product in $\Cal G$
and the standard product in $\Cal R[G]$.

The category structure shows up in the reassociating and (when $d>1$)
commuting isomorphisms. Specifically the
isomorphism $\alpha\:(a_f\sq a_g)\sq a_h\simeq a_f\sq( a_g\sq a_h)$ comes
from the thickening of the cone on four
points decomposed in two ways as union of cones on three points. The two
decompositions give two basepoints in
$[Y/\hat\partial Y,E]_0$, namely the homotopies $m_{fg,h}m_{f,g}$ and
$m_{f,gh}m_{g,h}$. These differ by a unit in
$R$, which gives the difference between the identifications of the iterated
products with $a_{fgh}$. But according to
2.3.2 this unit is exactly the associativity isomorphism in the category
associated with $E$. Thus the natural
isomorphism between the products in $\Cal G$  and $\Cal R[G]$ takes the
associativity isomorphisms in $\Cal G$ to the
$E$--twisted ones in $\Cal R[G]$.

A similar argument shows that the commutativity isomorphisms agree too,
when $d>1$.
\qed\enddemo

\head Modular field theories on 3--manifolds\endhead
Modular theories on 3--manifolds with a little extra data can be obtained as
follows: start with a theory on
4--dimensional thickenings of 2--complexes, corresponding to some
braided--symmetric category. Restrict to a subcategory
of objects that are almost determined by their boundaries. Then normalize
using an Euler-characteristic theory to
remove most of the remaining dependence on interiors. Here we carry this
through for group--categories. The untwisted
theories (which are $H_1$ theories in the sense of Section 3) can be
normalized if the order of the underlying group
is invertible. For cyclic groups we determine exactly which
group--categories give normalizable theories: in most cases
it requires a certain divisor of the group order to have a square root.
However there are cases, including the
 category with group $Z/2Z$ and $\sigma=-1$, that cannot be normalized.

\subhead 5.1 Extended, or weighted, 3--manifolds\endsubhead
There is a domain category (in the sense of \cite{Q}) with
\roster\item corners are closed 1--manifolds, with a parametrization of each
component by $S^1$;
\item boundaries are oriented surfaces with boundary, the boundary is a
corner object (ie, has parameterized
components, with correct orientation), and a lagrangian subspace of
$H_1(\overline{Y};Z)$; and
\item spacetimes are 3--manifolds whose boundaries are boundary objects
(ie, have lagrangian subspaces), together with
an integer (the ``index'').
\endroster
In (2) $\overline{Y}$ denotes the closed surface obtained by glueing copies
of $D^2$ to $Y$ via the given parameterizations
of the boundary components. A ``lagrangian subspace'' is a $Z$--summand of
half the rank on which the intersection
pairing vanishes. These objects are the ``extended'' or ``e--manifolds'' of
Walker \cite{W}, and special cases of
the ``weighted'' manifolds of Turaev \cite{T}. Turaev allows lagrangian
subspaces of the real rather than integer
cohomology.

 A domain category comes with cylinder functors and glueing operations.
Most of these are pretty clear. For instance
when glueing spacetimes along closed (no corners) boundaries, the weights
add. Glueing when corners are involved
requires Wall's formula for modified additivity of the index, using the
Shale--Weyl cocycle [\Cite{W}, \Cite{T}].

The geometric basis for the construction is:
\proclaim{5.1.1 Theorem} \roster\item If $U$ is an oriented 3--dimensional
thickening of a 1--complex, then the kernel
of the inclusion $H_1(\partial U,Z)\to H_1(U;Z)$ is a lagrangian subspace.
Every lagrangian subspace arises this way,
and the manifold $U$ is unique up to diffeomorphism rel boundary; and
\item (Kirby \cite{K}) A connected oriented 3--manifold is the boundary of a
smooth
4--manifold with the homotopy type of a 1--point union of copies of $S^2$. If
$X_1$ and $X_2$ are two such manifolds
with the same boundary, then for some $m_1,n_1,m_2,n_2$ there is a
diffeomorphism
$$X_1\#m_1CP^2\#n_1\overline{CP}^2 \simeq X_2\#m_2CP^2\#n_2\overline{CP}^2 $$
which is the identity on the boundary.\endroster\endproclaim

Some of the modifications in (2) can be tracked with the index of the
4--manifold: adding $CP^2$ increases it by 1,
while $\overline{CP}^2$ decreases it by 1. Doing both leaves the index
unchanged. This gives a refinement of (2):
\proclaim{5.1.2 Corollary} Suppose $X_1$ and $X_2$ are 4--manifolds as in
5.1.1(2) and the indexes are the same. Then
for some $p_1, p_2$ there is a diffeomorphism
$$X_1\#p_1(CP^2\#\overline{CP}^2) \simeq X_2\#p_2(CP^2\#\overline{CP}^2). $$
\endproclaim

\subhead 5.2 Construction of field theories\endsubhead
Now suppose $Z$ is a field theory on 4--dimensional thickenings of
2--complexes. Suppose $Y$ is a extended boundary object, so a surface with
parameterized boundary and homology
lagrangian. According to 5.1.1(1) this data is the same as a 3-d thickening
$U$ of a 1--complex with $\partial
U=\overline{Y}$, together with parameterized 2--disks in the boundary. This
is a boundary object of the  category of
thickenings, so we can define $\hat Z(Y)=Z(U)$.

If induced homomorphisms
$Z_V$  are unchanged by connected sum with $CP^2$ and
$\overline{CP}^2$ then we can define $\hat Z_X$ to be $Z_V$ for one of the
4--manifolds of 5.1.1(2) with $\partial V=X$.
Usually these operations do change $Z_V$; specifically there are elements
$\tau, \overline\tau\in R$ so
that
$$Z_{V\#CP^2}=\tau Z_V\text{\quad and \quad}
Z_{V\#\overline{CP}^2}=\overline{\tau} Z_V.\tag{5.2.1}$$
 These changes were called ``anomalies'' by
physicists. Usually the changes are too strong to fix, but
 sometimes we can fix the changes caused by
adding both $CP^2$ and $\overline{CP}^2$ together. Specifically, suppose
there is an inverse square root for
$\tau\overline{\tau}$: an element $r$ such that
$$ r^2\tau\overline\tau = 1.\tag{5.2.2}$$
Connected sum with $CP^2\#\overline{CP}^2$ changes $Z_V$ by
$\tau\overline\tau$ and increases the Euler characteristic
of $V$ by 2. Thus if we multiply by $r$ to the power $\Cal X(V)$ the
changes cancel. More specifically if
$(X,n)\:Y_1\to Y_2$ is a bordism in the extended 3--manifold category, and
$V\:U_1\to U_2$ is a corresponding 4-d
morphism of thickenings with index $n$ define
$$\hat Z_X=r^{\Cal X(V,U_1)}Z_{V}.\tag{5.2.3}$$

\proclaim{Proposition} If an element $r$ satisfying 5.2.2 exists, then
$\hat Z$ is a modular field theory on extended
3--manifolds.\endproclaim
Note that adding 1 to the index of an extended 3--manifold $X$ corresponds
to changing the bounding 4--manifold by
$\#CP^2$. This adds 1 to the Euler characteristic so changes $\hat Z_X$ by
$r\tau$. This is the ``anomaly'' of the
normalized theory. In particular it is nontrivial if $\tau\neq
\overline{\tau}$.
\demo{Proof} Multiplication by $r$ raised to the relative Euler
chacteristic gives a modular field theory with all
state spaces $R$, defined  on all finite complexes \cite{Q}. The product in
5.2.3 is the tensor
product of this Euler theory with
$Z$, so defines a theory on 4-d thickenings. Restricting to the
simply-connected thickenings obtained from extended
3--manifolds therefore is a modular field theory. By construction it is
insensitive to the difference between different
$V$ with fixed boundary and index, so it is a well-defined theory on
extended 3--manifolds.\qed\enddemo

\subhead 5.3 Normalization of group--category field theories\endsubhead
Here we describe the ``anomalies'' of the field theory associated to a
group--category in terms of the category
structure. In the cyclic case this is explicit enough to completely
determine when the field theory can be normalized
to give one on extended 3--manifolds.
\proclaim{5.3.1 Proposition} Suppose $\Cal G$ is a braided--commutative
group--category over $R$, with finite
underlying group $G$.  Then the associated field theory $Z$ has
$Z_{CP^2}=\Sigma_{g\in G}\sigma_g$ and $Z_{\overline{CP}^2}=\Sigma_{g\in
G}\sigma_g^{-1}$. If $G$ is cyclic of order
$n$,
$\sigma$ ($=\sigma_g$ for some generator $g$) has order exactly $\ell$, and
$R$ has no zero divisors then
$$Z_{CP^2\#\overline{CP}^2}=\cases n^2/\ell&\text{\quad if $\ell$ is odd}\\
2n^2/\ell&\text{\quad if $4|\ell$}\\
0&\text{\quad otherwise ($\ell$ is even and $\ell/2$ is odd).}\endcases$$
\endproclaim
We recall  $\sigma_g\in\units(R)$ is the number so that the commuting
endomorphism $\sigma_{g,g}\:g\sq g\to g\sq g$ is
multiplication by $\sigma_g$. The order of $\sigma_g$ divides the order of
$g$ if this order is odd, and twice this
order if it is even.
\subsubhead 5.3.2 Example\endsubsubhead
 If  $G=Z/2Z$ then $Z_{CP^2}=1+\sigma$ and
$Z_{\overline{CP}^2}=1+\sigma^{-1}$. Since $\sigma^4=1$  there are three
cases: $\sigma=1$, $\sigma=-1$, and $\sigma =i$ (a primitive $4^{th}$ root
of unity).
\roster
\item When $\sigma=1$ (the standard untwisted category) both $Z$ are 2, so
the inverse square root of the product is $1/2$. Thus the theory can be
normalized over $R[1/2]$ and gives an
anomaly-free theory ($\hat Z_X$ doesn't depend on the index of $X$).
\item When $\sigma=-1$ (the nontrivial symmetric category) both $Z$ are 0,
and no extended 3--manifold theory can be
obtained.
\item When $\sigma =i$ (a non-symmetric braided category) the $Z$ are $1+i$
and $1-i$ respectively. The product is
2, so the theory can be normalized over $R[1/\sqrt 2]$.
\endroster
Note that the  $Z/2Z$ category with $\sigma=-1$ is a (possibly twisted)
tensor factor of the quantum categories coming
from
$sl(2)$ at roots of unity. This should mean that on 4-d thickenings the
field theory is a (possibly twisted) tensor
product. The non-normalizability of the $Z/2Z$ factor would explain why it
has been so hard to normalize the full
$sl(2)$ theory.

\demo{Proof of 5.3.1} In general we want $Z_{CP^2-D^4}$, where $CP^2-D^4$
is regarded as a
bordism  $D^3\to D^3$ (relative to the corner $S^2=\partial D^3$). In the
group--category case this is the same as the
closed case ($CP^2$ as a bordism from the empty set to itself). This can
either be seen directly, or more generally
induced homomorphisms can be seen to be multiplicative with respect to
connected sums. Thus we consider the closed
case.

  Let
$k\in H^4(B^2_G;\units(R))$ be the class corresponding to the
group--category $\Cal G$. $Z_{CP^2}$
is multiplication by the sum over $[CP^2,B_G^2]=H^2(CP^2;G)=G$ of  $k$
evaluated on the image of the
fundamental class of
$CP^2$. We claim this evaluation for a single $g\in G$ is $\sigma_g$, so
the sum is as indicated in 5.3.1. The element
for $\overline{CP}^2$ is obtained by evaluating on the negative of the
fundamental class of $CP^2$, so gives
$\sigma_g^{-1}$.

This claim is verified using a geometric argument and the description of
2.3.4. Suppose data
$\tilde g\:D^2/S^1\to E$ and $\tilde m_{g,h}$ has been chosen. Then
$\sigma_{g,g}$ is obtained by glueing together
$\tilde m_{g,g}$, its inverse, and the standard commuting homotopy in
$\pi_2$ to get $D^2\times S^1\to E$. Consider
this as a neighborhood of a standard circle in $D^3$ and extend to
$D^3/S^2\to E$ by taking the complement to the
basepoint. $\sigma_{g,g}$ is the resulting element in $\pi_3(E)=\units(R)$.
We manipulate this a little.  $\tilde
m_{g,g}$ and its inverse cancel to leave just the standard commuting
homotopy. This gives the following description:
take an embedding
$\mu\:D^2\times S^1\to D^2\times S^1$ that goes twice around the
$S^1$, and locally preserves products. The element of $\pi_3$ is obtained by
$$D^3/S^2@>>> D^2\times S^1/(\partial D^2\times S^1)@> \mu^{-1}>>D^2\times
S^1/(\partial D^2\times
S^1)@>p>>D^2/S^1@>\tilde g>>E$$ where the first map divides out the
complement of  the standard $D^2\times S^1\subset
D^3$, and $p$ projects to the $D^2$ factor of the product. According to
2.3.4 the image of this in
$\pi_3(E)=\units(R)$ is $\sigma_g$.  Denote the composition
$D^3/S^2\to D^2/S^1$ by $h$. This is homotopic to the Hopf map. This can be
checked using Hopf's original
description: the inverses of two points in the interior of
$D^2$ give two unknotted circles in $S^3$ with linking number  1. But
$S^2\cup_h D^4\simeq CP^2$.
The vanishing higher homotopy of $B^2_G$ implies there is an extension
(unique up to homotopy) of $g$ over
the 4--cell to give $CP^2\to B^2_G$. General principles imply that the
$k$--invariant evaluated on the homology image of
the 4--cell (the orientation class of $CP^2$) is equal to the homotopy class
of the attaching map in $\pi_3E$, so the
evaluation does give $\sigma_g$.

The numerical presentation material of 2.5 can be used to make these
conclusions more concrete. We carry this out for
cyclic groups. Suppose $G$ is cyclic of order $n$ with generator $g$, and
$\sigma=\sigma_g$. Suppose
$\sigma$ has order $\ell$. From 2.5 we know $\ell$ divides $n$ if $n$ is
odd, and $2n$ if $n$ is even. Further (see
2.5.2), $\sigma_{g^r}=(\sigma_g)^{r^2}$. Therefore
$$Z_{CP^2}=\Sigma_{r=0}^{n-1}\sigma^{r^2}\text{\quad and
\quad}Z_{\overline{CP}^2}=\Sigma_{r=0}^{n-1}\sigma^{-r^2}.$$
The product of these is
$$\Sigma_{r,s=0}^{n-1}\sigma^{r^2-s^2}=\sigma_{r,s=0}^{n-1}\sigma^{(r+s)(r-s)}$$
Reindex this by setting $r-s= t$, and use the fact that $r$ and $s$ can be
changed by multiples of $n$ to get
$$\Sigma_{s,t=0}^{n-1}\sigma^{t(t+2s)}
=\Sigma_{t=0}^{n-1}\sigma^{t^2}\Sigma_{s=0}^{n-1}(\sigma^{2t})^s.\tag{5.3.3}$$ 
We have assumed $R$ has no zero divisors. This means if
$(\rho^n-1)=(\rho-1)(\Sigma_{s=0}^{n-1}\rho^s)=0$ then one of
the factors is 0. This implies
$$\Sigma_{s=0}^{n-1}\rho^s)=\cases n&\text{\quad if }\rho=1\\
0&\text{\quad if }\rho\neq 1.\endcases$$
Applying this with $\rho=\sigma^{2t}$ to (5.3.3) gives the sum over $t$
with $\sigma^{2t}=1$ of $n\sigma^{t^2}$.

If
$\ell$ (the order of $\sigma$) is odd or divisible by 4,  then
$\sigma^{2t}=1$ implies
$\sigma^{t^2}=1$. In this case the sum is just $n$ times the number of such
$t$ between 0 and $n-1$. This number is
$n/\ell$ if
$n$ is odd,
$2n/\ell$ if
$n$ is even. This gives the conclusion of the proposition in these cases.
If $\ell$ is even but $\ell/2$ is odd then
$\sigma^{2t}=1$ implies $\sigma^{t^2}=1$ if $t$ is even, and
$\sigma^{t^2}=-1$ if $t$ is odd. The sum is thus $n$
times the difference between the number of even and odd $t$ with
$\sigma^{2t}=1$. These are $t=(\ell/2)j$ for $0\leq
j<2n/\ell$, so they exactly cancel if $2n/\ell$ is even, or equivalently if
$\ell$ divides $n$. We are in the case
with $\ell$ even so $n$ is even and 4 divides $2n$. But $\ell/2$ odd, so if
$\ell$ divides $2n$ it must also divide
$n$. This completes the proof of the proposition.
\qed\enddemo

\Addresses\recd


\Refs


\ref\key {BFSV}\by C Balteany\by Z Fiedorowicz\by R Schw\"anzl\by R Vogt
\paper Iterated monoidal categories\paperinfo
preprint (1998)\endref

\ref\key {B}\by Ivelina Bobtcheva\paper On Quinn's invariants of
2--dimensional CW complexes\paperinfo preprint (1998)\endref

\ref\key {BQ}\by I Bobtcheva\by F Quinn\paper Numerical presentations of
tortile categories \paperinfo preprint (1997)\endref

\ref\key {CFS}\by J\,S Carter\by D\,E Flath\by M Saito\book The classical
and quantum $6i$--symbols\publ Princeton
University Press,\bookinfo Math.  Notes 43\yr 1996\endref

\ref\key {DW}\by R Dijkgraaf\by E Witten\paper Topological gauge theories
and group cohomology\jour Comm. Math. Physics,\vol129 \yr 1990\pages
393--429\endref

\ref\key {F}\by Z Fiedorowicz\paper The symmetric bar construction\paperinfo
preprint\endref

\ref\key {FQ}\by D Freed\by F Quinn\paper Chern--Simons theory with finite
gauge group\jour Comm. Math. Physics,\vol 156\yr 1993\pages435--472\endref

\ref\key {FK}\by J Fr\"olich\by T Kerler\book Quantum groups, quantum
categories, and quantum field
theories\bookinfo Lecture Notes in Mathematics,\vol 1542,\publ
Springer--Verlag\yr 1993\endref

\ref\key {GK}\by S Gelfand\by D Kazhdan\paper Examples of tensor
categories\jour Invent. Math. \vol 109\yr 1992\pages 595--617\endref


\ref\key {Ka}\by Christian Kassel\book Quantum Groups\bookinfo
Springer--Verlag Graduate Text in Mathematics,\vol 155\yr
1995\endref

\ref\key {Ke}\by T Kerler\paper Genealogy of non-perturbative
quantum-invariants of $3$--manifolds: the surgical
family\inbook Geometry and physics (Aarhus, 1995)\bookinfo
Lecture Notes in Pure and Appl. Math.\vol
184,\publ Dekker\yr 1997\pages 503--547\endref

\ref\key {K}\by R Kirby\paper A calculus for framed links in $S^3$\jour
Invent. Math. \vol45\yr 1978\pages 35--56\endref

\ref\key {Ku}\by G Kuperberg\paper Spiders for rank $2$ Lie algebras \jour
Comm. Math. Phys. \vol 180
\yr 1996\pages 109--151\endref

\ref\key {L}\by G Lusztig\book Introduction to quantum groups\bookinfo
Progress in Math. \vol 110,\publ Birkh\"auser\yr
1993\endref


\ref\key {McL}\by Saunders MacLane\book \book Categories for the working
mathematician\bookinfo Graduate texts in mathematics,
\vol 5,\publ Springer--Verlag\yr 1971\endref

\ref\key {M}\by S Majid\paper Representations, duals and quantum doubles of
monoidal categories\jour Rend. Circ. Mat.
Palermo, Suppl.
\vol26
\yr 1991\pages 197--206
\endref

\ref\key {Q}\by Frank Quinn\paper Lectures on axiomatic topological quantum
field theory\inbook IAS/Park City Math
Series\vol 1,\publ Amer. Math Soc.\yr 1995\pages 325--453\endref

\ref\key {QTP}\by Quantum Topology Project\publ{\tt
http://www.math.vt.edu/quantum\_topology} \endref

\ref\key {RT}\by N Yu Reshetikhin\by V\,G Turaev\paper Invariants of
3--manifolds via link polynomials and quantum groups\jour Invent. Math. \vol
103\yr1991\pages547--598\endref

\ref\key {Sa}\by Stephen Sawin\paper Links, quantum groups and TQFTs\jour
Bull. Amer. Math Soc. \vol 33\yr 1996\pages
413--445
\endref

\ref\key {S}\by M\,C Shum\paper Tortile tensor categories\jour J. Pure Appl.
Algebra, \vol93 \yr1994\pages 57--110
\endref

\ref\key {T}\by V Turaev\book Quantum invariants of knots and
3--manifolds\publ W de Gruyter,\publaddr Berlin\yr
1994\endref

\ref\key {W}\by K Walker\paper On Witten's 3--manifold invariants\paperinfo
preprint draft\yr 1991\endref

\ref\key {Wh}\by G\,W Whitehead\book Elements of homotopy theory\bookinfo
Graduate texts in Mathematics,\vol 61,\publ Springer--Verlag\yr
1978\endref

\ref \key {Y}\by D\,N Yetter\paper  Topological quantum field theories
associated to finite groups and crossed
$G$--sets\jour J. Knot Theory and its Ram.\vol 1 \yr 1992\pages 1--20\endref



\endRefs
\bye